\documentclass[a4paper,12pt,reqno]{amsart}
\usepackage{color}
\usepackage{txfonts}
\usepackage{bbm}
\usepackage{cases}
\usepackage{graphicx}
\usepackage{mathrsfs}
\usepackage{ulem}
\usepackage{stmaryrd}
\usepackage{amsfonts}
\usepackage{enumerate,amsmath,amssymb,amsthm,booktabs,url}
\usepackage{mathrsfs}
\usepackage[nobysame,abbrev]{amsrefs}
\BibSpec{article}{%
+{}{\PrintAuthors} {author}
+{,}{ \textrm} {title}
+{.}{ \textit} {journal}
+{,}{ \textbf} {volume}
+{}{ \parenthesize} {date}
+{,}{ } {pages}
+{.}{ Available at arXiv:} {eprint}
+{.}{} {transition}
}
\BibSpec{book}{%
+{}{\PrintAuthors} {author}
+{,}{ \textit} {title}
+{,}{ \textit} {edition}
+{.}{ \textrm} {series} %+{,}{ \textrm} {series}
+{,}{ Vol.} {volume} 
+{.}{ } {publisher}
+{,}{ } {date}
+{.}{} {transition}
}
\usepackage{color}
\usepackage[colorlinks,linktocpage,hyperindex]{hyperref}
\hypersetup{
    linkcolor=blue,          % color of internal links
    citecolor=red,        % color of links to bibliography
    filecolor=blue,      % color of file links
    urlcolor=cyan  
}
\newcommand{\blue}{\color{blue}}

\numberwithin{equation}{section}
\usepackage{thmtools} 
\declaretheorem[name=Theorem, parent=section]{theorem}
\declaretheorem[name=Lemma, sibling=theorem]{lemma}
\declaretheorem[name=Definition, sibling=theorem]{definition}

\declaretheorem[name=Corollary, sibling=theorem]{corollary}
\declaretheorem[name=Remark, sibling=theorem]{remark}
\declaretheorem[name=Example,sibling=theorem]{example}\declaretheorem[name=Proposition,sibling=theorem]{proposition}
\def\bt{\begin{theorem}}
\def\et{\end{theorem}}
\def\bl{\begin{lemma}}
\def\el{\end{lemma}}
\def\bd{\begin{definition}}
\def\ed{\end{definition}}
\def\bp{\begin{proposition}}
\def\ep{\end{proposition}}
\def\bc{\begin{corollary}}
\def\ec{\end{corollary}}
\def\br{\begin{remark}}
\def\er{\end{remark}}
\def\bexa{\begin{example}}
\def\eexa{\end{example}}
\def\bB{{\mathbf B}}
\def\bC{{\mathbf C}}

\def\1{{\mathbf{1}}}

\def\cA{{\mathcal A}}
\def\cB{{\mathcal B}}

\def\cF{{\mathcal F}}

\def\cR{{\mathcal R}}

\def\mB{{\mathbb B}}

\def\mD{{\mathbb D}}
\def\mE{{\mathbb E}}

\def\mL{{\mathbb L}}

\def\mN{{\mathbb N}}

\def\mP{{\mathbb P}}
\def\mQ{{\mathbb Q}}
\def\mR{{\mathbb R}}
\def\mS{{\mathbb S}}

\def\sB{{\mathscr B}}

\def\sF{{\mathscr F}}

\def\sI{{\mathscr I}}
\def\sJ{{\mathscr J}}

\def\sS{{\mathscr S}}

\def\p{\partial}
\def\geq{\geqslant}
\def\leq{\leqslant}
\def\ge{\geqslant}
\def\le{\leqslant}
\def\eps{\varepsilon}

\def\[{{\Big[}}
\def\]{{\Big]}}
\def\<{{\langle}}
\def\>{{\rangle}}
\def\({{\Big(}}
\def\){{\Big)}}

\def\e{{\rm e}}

\def\0{{\mathbf{0}}}

\def\dif{{\mathord{{\rm d}}}}

\def\div{\mathord{{\rm div}}}

\def\={&\!\!=\!\!&}

\usepackage{titletoc}
\titlecontents{section}[0pt]{\addvspace{1pt}\filright}
{{\thecontentslabel\quad}}
{ }{\titlerule*[6pt]{.}\contentspage}

\titlecontents{subsection}[1.7em]{\addvspace{1pt}\filright}
{{\thecontentslabel\quad}}
{ }{\titlerule*[6pt]{.}\contentspage}

\titlecontents{subsubsection}[3.4em]{\addvspace{1pt}\filright}
{{\thecontentslabel\quad}}
{ }{\titlerule*[6pt]{.}\contentspage}
\usepackage{titlesec}
\titleformat{\section}{\sc\filcenter}{\thesection.}{0.5em}{}[]
\titleformat{\subsection}{}{\thesubsection.}{0.5em}{\bfseries}[]
\titleformat{\subsubsection}{\it}{\thesubsubsection.}{0.5em}{}[]
\allowdisplaybreaks

\def\bpf{\begin{proof}}
\def\epf{\end{proof}}

\allowdisplaybreaks
\begin{document}

\title{Well-posedness of density dependent SDE driven by $\alpha$-stable process with H\"{o}lder drifts}

\date{\today}
\author{Mingyan Wu and Zimo Hao*}


\address{
Mingyan Wu: School of Mathematics and Statistics, Huazhong University of Science and Technology, Wuhan, Hubei 430072, P. R. China, Email: mingyanwu@hust.edu.cn}

\address{
Zimo Hao: School of Mathematics and Statistics, Wuhan University, Wuhan, Hubei 430072, P. R. China; 
Fakult\"at f\"ur Mathematik, Universit\"at Bielefeld, 33615, Bielefeld, Germany,
Email: zimohao@whu.edu.cn}

\thanks{Mingyan Wu is partially supported by National Natural Science Foundation of China (Grant No. 61873320). 
Zimo Hao is grateful to the financial supports of NNSFC grants of China (Nos. 12131019, 11731009), and the DFG through the CRC 1283 
``Taming uncertainty and profiting from randomness and low regularity in analysis, stochastics and their applications''. }
\thanks{* Corresponding author}

\begin{abstract}
In this paper, we show the weak and strong well-posedness of density dependent stochastic differential equations driven by  $\alpha$-stable processes with $\alpha\in(1,2)$. The existence part is based on Euler's approximation as  \cite{HRZ20}, while, the uniqueness is based on the Schauder estimates in Besov spaces for nonlocal Fokker-Planck equations. For the existence, we only assume the drift being continuous in the density variable.
For the weak uniqueness, the drift is assumed to be Lipschitz in the density variable, while for the strong uniqueness, we also need to assume the drift being 
$\beta_0$-order H\"older continuous in the spatial variable, where $\beta_0\in(1-\alpha/2,1)$.
\end{abstract}

\maketitle

%\tableofcontents

\section{Introduction}

Fix $\alpha \in (1,2)$. Let $(L_t)_{t\geq 0}$ be a  $d$-dimensional symmetric and rotationally invariant $\alpha$-stable process on some probability space $(\Omega,\sF,\mP)$. In this paper, we consider the following density dependent stochastic differential equation (abbreviated as DDSDE): 
\begin{align}\label{eq:AB}
\dif X_t = b(t,X_t,\rho_t(X_t)) \dif t + \dif L_t,\ \ X_0 \stackrel{(d)}{=} \mu_0,
\end{align}
where $b : \mR_+ \times \mR^d \times \mR_+ \to \mR^d$ is a bounded Borel measurable vector field,  $\mu_0$ is a probability measure over $\mR^d$ and for $t>0$, $\rho_t(x) = \mP\circ X_t^{-1}(\dif x )/\dif x$ is the distributional density of $X_t$ with respect to the Lebesgue measure $\dif x$ on $\mR^d$.

\medskip

In literature, DDSDE \eqref{eq:AB} is also called McKean-Vlasov SDE of Nemytskii-type which was firstly introduced in \cite[Section 2]{BR18} to give a probabilistic representation for the solutions of nonlinear Fokker-Planck equations. In a series of works \cites{BR18,BR20,BR21a,BR21b}, Barbu and R\"ockner investigated the following DDSDE driven by Brownian motions:
\begin{align}\label{eq:HG01}
\dif X_t=b(t,X_t,\rho_t(X_t))\dif t+\sigma(t,X_t,\rho_t(X_t))\dif W_t,\ \ X_0\stackrel{(d)}{=} \mu_0,
\end{align}
where $\sigma:\mR_+\times \mR^d \times \mR_+ \to \mR^{d}\otimes \mR^d$ is  measurable and $W$ is a standard $d$-dimensional Brownian motion. By It\^o's formula, one sees that $\rho_t$ solves the following nonlinear Fokker-Planck equation (NFPE) in the distributional sense:
 \begin{align*}
\p_t\rho_t-\frac{1}{2}\sum_{i,j=1}^d\p_i\p_j\[a_{ij}(t,\cdot,\rho_t)\rho_t\]+\div(b(t,\cdot,\rho_t)\rho_t)=0, \ \ \lim_{t\to0}\rho_t(x)\dif x=\mu_0(\dif x) \ \ \hbox{weakly},
\end{align*}
where $\p_i := \frac{\p}{\p_{x_i}}$, $a_{ij}:=\sum_{k=1}^d\sigma_{ik}\sigma_{jk}$, and $\div$ stands for the divergence.
More precisely, for any $\varphi\in C^\infty_0(\mR^d)$,
\begin{align*}
\<\rho_t,\varphi\>=\<\mu_0,\varphi\>+\frac{1}{2}\sum_{i,j=1}^d\int_0^t\<\rho_s,{ a_{ij}(s,\cdot,\rho_s) \p_{i}\p_{j}}\varphi\>\dif s+\int_0^t\<\rho_s,b(s,\cdot,\rho_s)\cdot\nabla \varphi\>\dif s,
\end{align*}
where 
$$
\<\rho_t,\varphi\>:=\int_{\mR^d}\varphi(x)\rho_t(x)\dif x=\mE\varphi(X_t).
$$
 In Barbu and R\"ockner's works, they obtained the well-posedness for NFPE through analytic methods, and then used the so-called superposition principle to get the well-posedness of DDSDE \eqref{eq:HG01}. Recently, different from these works, the second named author together with R\"ockner and Zhang \cite{HRZ20} gave a purely probabilistic proof for the existence of the solution to 
 the following DDSDE with additive noises:
\begin{align}\label{eq:DG01}
\dif X_t=b(t,X_t,\rho_t(X_t))\dif t+\dif W_t.
\end{align}

It is well known that Brownian motion is a continuous L\'evy process.  Hence, it is natural to consider such density dependent SDEs driven by  pure jump L\'evy processes. In particular, we consider $\alpha$-stable processes which are typical L\'evy processes having selfsimilar properties (cf. \cite{Sa99}). Up to now, the study of the well-posedness of SDEs with stable noises has been and remains an important area in stochastic analysis. For the classical case, there are a lot of results about strong solutions, weak solutions, and martingale solutions (see \cite{Pr12}, \cite{LZ19}, \cite{CHZ20}, \cite{HWZ20} and etc.).  We also see that there are  many results about McKean-Vlasov SDEs with jumps (see \cite{LMW21} and references therein). Among these results, some applications can be found in financial mathematics (cf. \cite{BCD17}) and neural net-works (cf. \cite{Ma07}). However, under the framework of L\'evy noises, there is no any results about Nemytskii's type SDEs. Thus, it is natural and interesting to investigate DDSDE \eqref{eq:AB}.

\medskip

On the other hand, McKean-Vlasov SDEs with L\'evy noises { are} related to non-local integral-PDEs. By It\^o's formula (cf. \cite[Theorem 5.1]{IW89}) for DDSDE \eqref{eq:AB}, we have that for any $\varphi\in C^\infty_0(\mR^d)$,
\begin{align}\label{FPENL}
\<\rho_t,\varphi \> = \<\mu_0,\varphi \> + \int_0^t  \< \rho_s, b(s,\cdot,\rho_s)\cdot \nabla \varphi\>\dif s+ \int_0^t \<\rho_s,\Delta^{\alpha/2} \varphi \> \dif s,
\end{align}
where
\begin{align}\label{eq:AA06}
\begin{split}
\Delta^{\alpha/2} \varphi (x) 
& :=  \int_{\mR^d} \( \varphi(x+z) - \varphi(x) - z \1_{|z|\leq 1} \cdot \nabla \varphi(x) \) {|z|^{-d-\alpha}} \dif z\\
& = \frac{1}{2} \int_{\mR^d} \( \varphi(x+z) + \varphi(x-z) - 2\varphi(x) \) {|z|^{-d-\alpha}} \dif z
\end{split}
\end{align}
is the infinitesimal  generator of $({ L_t})_{t \geq 0}$ (cf. \cite[Theorem 31.5]{Sa99}). Consequently, one sees that $\rho_t$ solves the following equation in the distributional sense:
\begin{align}\label{eq:FG01}
\p_t \rho_t - \Delta^{\alpha/2} \rho_t + \div (b(t,\cdot,\rho_t)\rho_t) = 0,\ \ \lim_{t \downarrow 0} \rho_t(x) \dif x= \mu_0(\dif x) \ \ \hbox{weakly},
\end{align}
where we use the fact that $ \Delta^{\alpha/2}$ is a self-adjoint operator. We point out that the infintesimal generator of Brownian motion is the Laplacian $\Delta$. The fractional Laplacian operator $\Delta^{\alpha/2}$ is non-local, and is essentially different from
the local operator $\Delta$. For instance, we can use Leibniz's rule to handle $\Delta(fg)$ but  the non-local case is more difficult. Thus, the Euler's type approximation in \cite{HRZ20}, a purely probabilistic method, is chosen to show the existence of the solutions of DDSDE \eqref{eq:AB} in this paper. 

\medskip

Moreover, when $b(t,\cdot,u)$ is $\beta_0$-order H\"older continuous uniformly in $t,u$ with $\beta_0\in(1-\alpha/2,1)$, we obatin the uniqueness based on some priori estimates of  Besov-type  (see Lemma \ref{lem:PK01}) for the nonlcal Fokker-Planck equation \eqref{FPENL}. This part is not studied in \cite{HRZ20}.  It is worth noting that the condition $\beta_0>1-\alpha/2$ is natural. The uniqueness in \cite{HRZ20} is obtained based on the well-known pathwise uniqueness for SDE \eqref{eq:DG01} with bounded measurable drift $b(t,x,\rho_t(x))$ (cf. \cite{Ve79}). However, the situtation changes when we consider $\alpha$-stable noises with $\alpha \in (0,2)$. Let us consider
$$
\dif X_t = b(t,X_t) \dif t + \dif L_t,
$$
where $L$ is a $d$-dimensional symmetric $\alpha$-stable process. When $d=1$ and $\alpha<1$, even a bounded and $\beta_0$-H\"older continuous $b$ is not enough to ensure pathwise uniqueness if $\alpha + \beta_0<1$ (see \cite{TTW} for the counterexample). When $d\geq 1$ and $\alpha\in[1,2)$, Priola \cite{Pr12} obtained the pathwise uniqueness under $\beta_0>1-\alpha/2$. The condition $\beta_0>1-\alpha/2$ can be found in \cite{CZZ21} and \cite{HWZ20} as well for the supercritical case and the degenerate case respectively.

\medskip

Before stating the main result, we introduce the classical H\"older spaces in $\mR^d$. For $\beta>0$, let $\bC^\beta(\mR^d)$ be the classical $\beta$-order H\"older space consisting of all measurable functions $f:\mR^d\to\mR$ with
\begin{align*}
\|f\|_{\bC^\beta}:=\sum_{j=0}^{[\beta]}\|\nabla^jf\|_\infty+[\nabla^{[\beta]}f]_{\bC^{\beta-[\beta]}}<\infty,
\end{align*}
where $[\beta]$ denotes the greatest integer less than $\beta$, $\nabla^j$ stands for the $j$-order gradient, and
\begin{align*}
\|f\|_\infty:=\sup_{x\in\mR^d}|f(x)|,\quad [f]_{\bC^\gamma}:=\sup_{h\in\mR^d}\frac{\|f(\cdot+h)-f(\cdot)\|_\infty}{|h|^\gamma},~\gamma\in(0,1).
\end{align*}
%For any $p\in [1,\infty]$, we use $\ell^p$ to denote the usual space of a sequence of real numbers that is $p$-order summable.
In the sequel, for any $p\in [1,\infty)$, we denote by $L^p$ the space of all $p$-order integrable functions on $\mR^d$ with the norm denoted by $\|\cdot\|_p$.%$\|\cdot \|_{\ell^p}$ 

\medskip

As mentioned before, to show the existence of a weak solution, we consider the following Euler scheme to DDSDE \eqref{eq:AB}: Let $T>0$, $N \in \mN$ and $h:=T/N$. For $t \in [0,h]$, define
$$
X_t^N := X_0 + L_t,
$$
and for $t \in (kh,(k+1)T]$ with $k=1,\cdots,N-1$, we inductively define $X_t^N$ by 
\begin{align*}
X_t^N := X_{kh}^N  + \int_{kh}^t b(s, X_{kh}^N, \rho_{kh}^N(X_{kh}^N)) \dif s + ( L_t - L_{kh} ),
\end{align*}
where $\rho_{kh}^N(x)$ is the distributional density of $X_{kh}^N$, whose existence is easily seen from the construction. 

\medskip

We give the definition of a weak solution to DDSDE  \eqref{eq:AB}:
 
\bd[Weak solutions]\label{In:Def}
Let $\mu_0$ be a probability measure on $\mR^d$ and $\alpha \in(1,2)$. We call a filtered probability space $(\Omega, \sF, \mP;( \sF_t)_{t\geq 0})$ together with a pair of $\sF_t$-adapted processes $( X_t, L_t)_{t \geq 0}$ defined on it a weak solution of SDE \eqref{eq:AB} with initial distribution $\mu_0$, if 
\begin{enumerate}[(i)]
\item $ \mP \circ X_0^{-1} = \mu_0$, and $ (L_t)_{t\geq 0}$ is a $d$-dimensional  { symmetric and } rotationally invariant $\alpha$-stable process.
\item for each $t>0$, $\rho_t(x)= \mP \circ X_t^{-1}(\dif x)/\dif x $ and
$$
 X_t = X_0+\int_0^t b(s, X_s,\rho_s( X_s)) \dif s +  L_t,\ \ \mP-{\rm a.s}.
$$
\end{enumerate}
\ed

The following existence and uniqueness result is the main theorem of this paper.

\bt\label{in:Main}
Assume that $\alpha\in(1,2)$, and $b$ is bounded measurable and for any $(t,x,u_0)\in\mR_+\times\mR^d\times\mR_+$,
\begin{align}\label{eq:PO}
\lim_{u \to u_0} | b(t,x,u) - b(t,x,u_0) | = 0.
\end{align}
{\bf (Existence)}
For any $T>0$ and initial distribution $\mu_0$, there are a subsequence $N_k$ and a weak solution $(X, L)$ to DDSDE \eqref{eq:AB} in the sense of Definition \ref{In:Def} so that for any bounded measurable $f$ on $\mR^d$ and $t\in(0,T]$,
\begin{align*}
\lim_{k\to\infty}\mE f(X^{N_k}_t)=\mE f(X_t).
\end{align*}
Moreover, for each $t\in(0,T]$, $X_t$ admits a density $\rho_t$ satisfying the estimate
\begin{align}\label{In:M1}
\rho_t(y)\le c \int_{\mR^d}\frac{t}{(t^{1/\alpha}+|x-y|)^{d+\alpha}}\mu_0(\dif x),
\end{align}
where the constant $c >0$ only depends on $T,d,\alpha,\|b\|_\infty$, and the following $L^1$-convergences hold:
\begin{align}\label{In:M2}
\lim_{k\to\infty}\int_{\mR^d}|\rho^{N_k}_t(y)-\rho_t(y)|\dif y=0
\end{align}
and
\begin{align}\label{In:M3}
\lim_{k\to\infty}\int_0^T\int_{\mR^d}|\rho^{N_k}_t(y)-\rho_t(y)|\dif y\dif t=0.
\end{align}
{\bf (Uniqueness)}
Suppose that there is a constant $c >0$ such that for all $(t,x,u_i)\in\mR_+\times\mR^d\times\mR_+$, $i=1,2$,
\begin{align}\label{Lipb}
| b(t,x,u_1)-b(t,x,u_2)|\leq c |u_1-u_2|.
\end{align}
{\bf i)} If $\mu_0(\dif x)=\rho_0(x)\dif x$ with $\rho_0\in L^q(\mR^d)$ for some $q\in(\frac{d}{\alpha-1},+\infty]$, then the weak uniqueness holds for DDSDE \eqref{eq:AB}.\\
{\bf ii)} { If $\mu_0(\dif x)=\rho_0(x)\dif x$ with $\rho_0\in \bC^{\beta_0}(\mR^d)$ for some $\beta_0\in(1-\alpha/2,1)$ and }
$$
\sup_{(t,u)\in\mR_+^2}\|b(t,\cdot,u)\|_{\bC^{\beta_0}}<\infty,
$$
then the strong uniqueness holds for DDSDE \eqref{eq:AB}. 
\et

%
%%
%%%
%%%%

\iffalse%%%%
\begin{proof}\blue
Notice that \eqref{In:M1} is obtained from \eqref{eq:RR01}, \eqref{eq:DD04-1}, \eqref{S2:Ep} and \eqref{eq:TT01}. Observing that, by  \eqref{Se3:Main}, \eqref{In:M1}, and \eqref{eq:DD04-1}, for any $M>0$,
\begin{align}\label{eq:ZH01}
\int_{\mR^d}|\rho^{N_k}_t(y)-\rho_t(y)|\dif y \leq \int_{|y|\leq M}|\rho^{N_k}_t(y)-\rho_t(y)|\dif y + 2\int_{|y|>M} \int_{\mR^d}t^{-d/\alpha} \mu_0(\dif x)\dif y.
\end{align}
Hence, we get \eqref{In:M2} by taking $k$ to $\infty$ and $M$ to $\infty$. Furthermore, for any $\eps>0$, 
 \begin{align*}
\sup_k\int_0^\eps\int_{\mR^d}|\rho^{N_k}_t(y)-\rho_t(y)|\dif y\dif t\le 2\eps. 
\end{align*}
Based on \eqref{eq:ZH01}, we obtain \eqref{In:M3} by taking $k$ to $\infty$ and $\eps\to0$. The proofs for existence and uniqueness are presented in Section \ref{Sec:Main}.
\end{proof}
\fi%%%%

%%%%
%%%
%%
%

\br
Although we use the same method as \cite{HRZ20} in the existence part, our assumptions on drifts are weaker. Therein, the following local uniform continuity is assumed,
\begin{align}\label{ZZ003}
\lim_{u\to u_0}\sup_{|x|<R}|b(t,x,u)-b(t,x,u_0)|=0,\ \ \forall  t\geq 0, R>0.
\end{align}
For example, $b(t,x,u):=(|u/x|\wedge1)\1_{x\ne0}$, $x\in\mR$, $u\in\mR_+$, satisfies  the condition \eqref{eq:PO} but doesn't satisfy \eqref{ZZ003} for any $R>0$.
\er

\br
For the uniqueness, the conditions here are natural. Conditions in \textbf{i)} are the same as \cite[Theorem 1.2]{HRZ20}; the same condition $\beta_0>1-\alpha/2$ in \textbf{ii)} can be found in \cites{CZZ21,Pr12} as well. 
\er

%\subsection{Overview of the paper and notations}

The paper is organized as follows. In Section \ref{Sec2}, we show some estimates of the density to the rotationally invariant and symmetric $\alpha$-stable process. In Section \ref{Sec:Be}, we introduce Besov spaces and establish Schauder's estimates for non-local parabolic equations by using Littlewood-Paley's type estimates of heat kernels. In Section \ref{Sec3}, we prove some uniform estimates in $N$ about heat kernels of Euler's scheme $X^N_t$. In Section \ref{Sec:Main}, we show the proof of Theorem \ref{in:Main}. 

\medskip

Throughout this paper, we use the following conventions and notations: As usual, we use $:=$ as a way of definition. Define $\mN_0:= \mN \cup \{0\}$ and $\mR_+:=[0,\infty)$. The letter $c=c(\cdots)$ denotes an unimportant constant, whose value may change in different palces. We use $A \asymp B$ and $A\lesssim B$ to denote $c^{-1} B \leq A \leq c B$ and $A \leq c B$, respectively, for some unimportant constant $c \geq 1$. 

\section{Preliminaries}\label{Sec2}

\subsection{$\alpha$-stable processes}

A c$\rm\grave{a}$dl$\rm\grave{a}$g process $\{L_t\mid t \geq 0\}$ on $\mR^d$ ($d\geq 1$) is called a L\'evy process, if $L_0 = 0$ almost surely and $L$ has independent and identically distributed increments. The associated Poisson random measure is defined by
$$
N((0,t]\times \Gamma)  :=  \sum_{s\in(0,t]} \1_{\Gamma}(L_s -  L_{s-}),\ \ \Gamma \in \sB(\mR^d\setminus \{0\}) , t >0,
$$
and the L\'evy measure is given by
\begin{align*}
\nu(\Gamma):=\mE N((0,1]\times\Gamma).
\end{align*}
Then,  the compensated Poisson random measure is defined by
$$
\tilde N(\dif r,\dif z): = N(\dif r,\dif z) - \nu(\dif z)\dif r.
$$ 
For $\alpha\in(0,2)$, a L\'evy process $L_t$ is called a symmetric and rotationally invariant $\alpha$-stable process if the L\'evy measure has the form
 $$
 \nu^{(\alpha)}(\dif z)= {c|z|^{-d-\alpha}} \dif z,
 $$
with some specific constant $c=c(d,\alpha)>0$. In this paper, we only cosider the symmetric and rotationally invariant $\alpha$-stable process. Without causing confusion,  we simply call it the $\alpha$-stable  process, and assume that $\nu^{(\alpha)}(\dif z)= {|z|^{-d-\alpha}} \dif z$ here and after. For  any  $0\le\gamma_1<\alpha< \gamma_2$, it is easy to see that
\begin{align}\label{eq:AA05}
\int_{\mR^d} (|z|^{\gamma_1} \wedge |z|^{\gamma_2}) \nu^{(\alpha)}(\dif z)< \infty.
\end{align}
By L\'evy-It\^o's decomposition (cf. \cite[Theorem 19.2]{Sa99}, we have
\begin{align}\label{eq:EE04}
L_t =  \lim_{\eps \downarrow 0} \int_0^{t} \int_{\eps<|z|\leq 1} z \tilde N(\dif r,\dif z) + \int_0^{t} \int_{|z|> 1} z   N(\dif r,\dif z).
\end{align}
By \cite[Theorem 31.5]{Sa99}, the infinitesimal  generator of L\'evy process $({ L_t})_{t \geq 0}$ is the fractional { Laplacian operator} $\Delta^{\alpha/2}$ defined by \eqref{eq:AA06}. 

\medskip

Moreover, by L\'evy-Khintchine's formula \cite[Theorem ]{Sa99}, for $\forall |\xi|\geq1$, we have 
\begin{align*}
|\mE e^{i\xi\cdot L_{t} }|\leq&\exp\left(t\int_{\mR^d}(\cos(\xi\cdot z)-1)\nu^{(\alpha)}(\dif z)\right)\\
\leq&\exp\left(-{t}|\xi|^\alpha\int_0^{\infty}\int_{\mS^{d-1}}\frac{1-\cos(\bar{\xi}\cdot r\theta)}{r^{1+\alpha}}\Sigma(\dif \theta)\dif r\right)
\leq e^{-ct|\xi|^{\alpha}},
\end{align*}
where $\Sigma$ is the uniform measure on the sphere $\mS^{d-1}$, and the constant $c>0$ depends only on $\alpha$ and $\Sigma(\mS^{d-1})$. Hence, by \cite[Proposition 28.1 ]{Sa99}, $L_t$ admits a smooth density function $p_\alpha(t,\cdot)$ given by Fourier's inverse transform
$$
p_\alpha(t,\cdot)=(2\pi)^{-d/2}\int_{\mR^d}e^{-i x\cdot \xi}\mE e^{i\xi\cdot L_{t} }\dif \xi,\ \ \forall t>0,
$$
and the partial derivatives of $p_\alpha(t,\cdot)$ at any orders  tend to $0$ as $|x|\to \infty$. Since the $\alpha$-stable  process $L_t$ has the scaling property
\begin{align*}
(\lambda^{-1/\alpha} L_{\lambda t})_{t\geq 0} \stackrel{(d)}{=} ( L_t )_{t \geq 0},\ \ \forall \lambda >0,
\end{align*}
it is easy to see that
\begin{align}\label{eq:HG01}
p_\alpha(t,x) = t^{-d/\alpha}p_\alpha(1,t^{-1/\alpha}x).
\end{align}
By \cite[Theorem 2.1]{BG60}, one knows that there is a constant $ c=c(d,\alpha)>1$ such that 
\begin{align}\label{eq:AA09} 
c^{-1}\varrho_\alpha(t,x)\leq p_\alpha(t,x) \leq c \, \varrho_\alpha(t,x), 
\end{align} 
where
\begin{align}\label{eq:TT01}
\varrho_\alpha(t,x) : = \frac{t}{(t^{1/\alpha} + |x|)^{d+\alpha}}.
\end{align}
By \cite{CZ16}*{Lemma 2.2}, for any $j\in\mN_0$, there is a constant $ c=c(j,\alpha)>0$ such that
\begin{align}\label{eq:AA08}
|\nabla^j p_\alpha(t,x)| \leq c t^{-j/\alpha}\varrho_\alpha(t,x).
\end{align}
Since
\begin{align*}%\label{eq:AA02}
(t^{1/\alpha} + |x+z|)^{-\gamma}\leq 4^\gamma (t^{1/\alpha} + |x|)^{-\gamma},\ \ \hbox{for}\ \  |z| \leq (2 t^{1/\alpha})\vee (|x|/2),
\end{align*}
we get that
\begin{align}\label{eq:AA02}
\varrho_\alpha (t,x+z) \leq 4^{d+\alpha} \varrho_\alpha (t,x),\ \ \hbox{for}\ \  |z| \leq (2 t^{1/\alpha})\vee (|x|/2).
\end{align}
Note that $p_\alpha(t,x)$ is the heat kernel of the operator $\Delta^{\alpha/2}$, i.e.,
\begin{align}\label{eq:AA11}
\p_t p_\alpha(t,x) = \Delta^{\alpha/2} p_\alpha(t,x),\ \ \lim_{t \downarrow 0} p_\alpha(t,x) = \delta_0(x),
\end{align}
where $\delta_0$ is the Dirac measure. We aslo have the following Chapman-Kolmogorov (abbreviated as C-K) equations:
\begin{align}\label{eq:CC00}
(p_\alpha(t)* p_\alpha(s))(x) = \int_{\mR^d} p_\alpha(t,x-y)p_\alpha(s,y) \dif y = p_\alpha(t+s,x),\ \ t,s>0.
\end{align}

\subsection{Some estimates of the heat kernel of $\Delta^{\alpha/2}$}

Now we give some estimates of the  heart kernel of $\Delta^{\frac{\alpha}{2}}$. These estimates are straightforward and elementary. Note that  Lemma \ref{lem:AA02} and Corollary \ref{cor:PQ01} are the same as \cite[Lemma 2.2]{CZ16} and \cite[Theorem 2.4]{CZ16} respectively when $j=0$.

\bl\label{lem:AA02}
For any $j \in \mN_0$ and $\beta \in (0,1)$, there is a constant $c= c(d,\alpha,\beta,j)>0$ such that for every $t>0,x_1,x_2\in\mR^d$,
\begin{align}\label{eq:AA01}
|\nabla^j p_\alpha(t,x_1) - \nabla^j p_\alpha(t,x_2)| \leq c |x_1-x_2|^\beta t^{-{(j+\beta)}/{\alpha}}(p_\alpha(t,x_1) + p_\alpha(t,x_2)).
\end{align}
\el
\begin{proof}
If $|x_1 - x_2| > t^{1/\alpha}$, then by \eqref{eq:AA08} we have
\begin{align*}
|\nabla^j p_\alpha(t,x_1) - \nabla^j p_\alpha(t,x_2)| 
&\lesssim t^{-j/\alpha}(\varrho_\alpha(t,x_1) + \varrho_\alpha(t,x_2))\\
&\lesssim |x_1-x_2|^\beta t^{-{(j+\beta)}/{\alpha}}(\varrho_\alpha(t,x_1) + \varrho_\alpha(t,x_2)).
\end{align*}
If $|x_1 - x_2| \leq t^{1/\alpha}$, then by the mean-value formula and \eqref{eq:AA08},
\begin{align*}
|\nabla^j p_\alpha(t,x_1) - \nabla^j p_\alpha(t,x_2)|  
& \leq |x_1-x_2|\int_0^1 |\nabla^{j+1} p_\alpha(t,x_1 + \theta(x_2-x_1)| \dif \theta\\
& \lesssim  |x_1-x_2 |t^{-(j+1)/\alpha} \int_0^1 \varrho_\alpha(t,x_1+ \theta(x_2-x_1))\dif \theta\\
&\stackrel{\eqref{eq:AA02}}{\lesssim} |x_1-x_2 |^\beta t^{-(j+\beta)/\alpha} \varrho_\alpha(t,x_1),
\end{align*}
where we have used $\beta \in (0,1)$ in the last inequality. Combining the above calculations, we get \eqref{eq:AA01} by \eqref{eq:AA09}.
\end{proof}

As a consequence of Lemma \ref{lem:AA02}, we have the following corollary. 

\bc\label{cor:PQ01}
For any $j \in \mN_0$, there is a constant $c= c(d,\alpha,j)>0$ such that for every $t>0$ and $ x\in\mR^d$,
\begin{align}\label{Iq:HK01}
|\Delta^{\frac{\alpha}{2}}\nabla^j p_\alpha(t,x)| \leq c t^{-1-j/\alpha} p_\alpha(t,x).
\end{align}
\ec
\bpf
First of all, recalling the definition \eqref{eq:AA06},
\begin{align*}
\Delta^{\frac{\alpha}{2}}\nabla^j p_\alpha(t,x)&=\int_{|z|\le t^{1/\alpha}}\(\nabla^j p_\alpha(t,x+z)-\nabla^j p_\alpha(t,x)-z\cdot\nabla^{j+1} p_\alpha(t,x)\)\frac{\dif z}{|z|^{d+\alpha}}\\
&+\int_{|z|\ge t^{1/\alpha}}\(\nabla^j p_\alpha(t,x+z)-\nabla^j p_\alpha(t,x)\)\frac{\dif z}{|z|^{d+\alpha}}\\
&:=\sI_1+\sI_2.
\end{align*}
For $\sI_1$, by \eqref{eq:AA01} and \eqref{eq:AA02}, we have that for any $\beta\in(\alpha-1,1)$,
\begin{align*}
\sI_1 \le \int_{|z|\le t^{1/\alpha}}|z|^{\beta+1}t^{-\frac{j+1+\beta}{\alpha}}\frac{\dif z}{|z|^{d+\alpha}}\varrho_\alpha(t,x)\lesssim t^{-1-j/\alpha}\varrho_\alpha(t,x).
\end{align*}
For $\sI_2$, by \eqref{eq:AA08}, ones see that
\begin{align*}
\sI_2&\lesssim   t^{-j/\alpha} \int_{|z|>t^{1/\alpha}}\(\varrho_\alpha(t,x+z)+\varrho_\alpha(t,x)\)\frac{\dif z}{|z|^{d+\alpha}}\\
&\lesssim t^{-1-j/\alpha}\varrho_\alpha(t,x)+t^{-j/\alpha} \int_{|z|>t^{1/\alpha}}\varrho_\alpha(t,x+z)\frac{\dif z}{|z|^{d+\alpha}}.
\end{align*}
Then, we only need to estimate the second term above denoted by $\sI_3$. If $|x|\le 2t^{1/\alpha}$, by \eqref{eq:AA02} and  \eqref{eq:TT01}, we obtain that
\begin{align*}
\sI_3 
&\lesssim t^{-j/\alpha} \int_{|z|>t^{1/\alpha}}\varrho_\alpha(t,z)\frac{\dif z}{|z|^{d+\alpha}}
\le t^{-j/\alpha}\frac{t}{t^{(d+\alpha)/\alpha}}\int_{|z|>t^{1/\alpha}}\frac{\dif z}{|z|^{d+\alpha}} \\
&\lesssim t^{-(d+\alpha+j)/\alpha}= t^{-j/\alpha-1}\frac{t}{t^{(d+\alpha)/\alpha}}
{\lesssim} t^{-j/\alpha-1}\varrho_\alpha(t,x).
\end{align*}
If $|x|> 2t^{1/\alpha}$, by \eqref{eq:AA02} and \eqref{eq:TT01}, we have
\begin{align*}
\sI_3&=t^{-j/\alpha} \(\int_{\frac{|x|}{2}\ge|z|>t^{1/\alpha}}+\int_{|z|>\frac{|x|}{2}}\)\varrho_\alpha(t,x+z)\frac{\dif z}{|z|^{d+\alpha}}\\
&\lesssim t^{-j/\alpha} \varrho_\alpha(t,x)\int_{|z|>t^{1/\alpha}}\frac{\dif z}{|z|^{d+\alpha}}+t^{-j/\alpha}\frac{1}{|x|^{d+\alpha}}\int_{|z|>\frac{|x|}{2}}\varrho_\alpha(t,x+z)\dif z\\
&\stackrel{\eqref{eq:AA09}}{\lesssim} t^{-j/\alpha-1}\varrho_\alpha(t,x) + t^{-j/\alpha}\frac{1}{|x|^{d+\alpha}} 
\lesssim t^{-j/\alpha-1}\varrho_\alpha(t,x).
\end{align*}
Based on \eqref{eq:AA09}, the proof is complete.
\epf

The following result is also ture when we consider Gaussian heat kernels (cf.  \cite[Lemma 2.1]{HRZ20}).

\bl\label{lem:AA03}
For any  $\beta \in (0,\alpha)$ and $j \in\mN_0$, there is a constant $c=c(d,\alpha,\beta,j)>0$ such that for every $t_1,t_2>0$ and $x\in \mR^d$,
\begin{align}\label{eq:AA03}
|\nabla^j p_\alpha(t_1,x) - \nabla^j p_\alpha(t_2,x)| 
\leq c |t_2-t_1|^{\beta/\alpha}   (  t_1^{-(j+\beta)/\alpha} p_\alpha(t_1,x) + t_2^{-(j+\beta)/\alpha} p_\alpha(t_2,x) ).
\end{align}
\el
\begin{proof}
Without loss of generality, we assume that $t_2>t_1$. If $t_2 - t_1> t_1$, then $t_1\vee t_2 \leq 2(t_2 -t_1)$ and 
\begin{align*}
|\nabla^j p_\alpha(t_1,x) - \nabla^j p_\alpha(t_2,x)|
& \stackrel{\eqref{eq:AA08}}{\lesssim} t_1^{-{j}/{\alpha}}\varrho_\alpha(t_1,x) + t_2^{-{j}/{\alpha}}\varrho_\alpha(t_2,x) \\
& \lesssim |t_2-t_1|^{\beta/\alpha}   (  t_1^{-(j+\beta)/\alpha}\varrho_\alpha(t_1,x) + t_2^{-(j+\beta)/\alpha}\varrho_\alpha(t_2,x) ).
\end{align*}
For $t_2 - t_1 \leq t_1$, notice that by \eqref{eq:AA11} and \eqref{Iq:HK01}, 
\begin{align}\label{eq:AA10}
|\nabla^j \p_t p_\alpha(t,x)| =  |\nabla^j  \Delta^{\alpha/2} p_\alpha(t,x)| = |\Delta^{\alpha/2} \nabla^j  p_\alpha(t,x)| \lesssim  t^{-1-j/\alpha}  \varrho_\alpha (t,x).
\end{align}
Thus, by the mean-value formula and $\beta \in (0,\alpha)$, we have
\begin{align*}
|\nabla^j p_\alpha(t_1,x) - \nabla^j p_\alpha(t_2,x)|
& \leq |t_2-t_1|\int_0^1 |\nabla_j \p_t p_\alpha|(t_1+\theta(t_2-t_1),x) \dif \theta\\
& \lesssim  |t_2-t_1| \int_0^1 (t_1+\theta(t_2-t_1))^{-1-j/\alpha}  \varrho_\alpha (t_1+\theta(t_2-t_1),x) \dif \theta\\
& \stackrel{\eqref{eq:TT01}}{\lesssim} |t_2-t_1| t_1^{-1-j/\alpha}  \varrho_\alpha (t_1,x)
\leq |t_2-t_1|^{\beta/\alpha} t_1^{-(j+\beta)/\alpha}  \varrho_\alpha (t_1,x). 
\end{align*}
By \eqref{eq:AA09}, the proof is finished.
\end{proof}

\section{Besov spaces and Schauder's estimates}\label{Sec:Be}

In this section, we introduce Besov spaces where we obtain Schauder's estimates for the operator $\p_t-\Delta^{\alpha/2}$ (see Lemma \ref{lem:PK01} below). Let $\sS(\mR^d)$ be the Schwartz space of all rapidly decreasing functions on $\mR^d$, and $\sS'(\mR^d)$ 
the dual space of $\sS(\mR^d)$ called Schwartz generalized function (or tempered distribution) space. Given $f\in\sS(\mR^d)$, 
the Fourier and inverse transforms of $f$ are defined by
$$
\hat f(\xi):=\cF f(\xi):=(2 \pi)^{-d/2}\int_{\mR^d} \e^{-i\xi\cdot x}f(x)\dif x, \quad\xi\in\mR^d
$$
and
$$
\check f(x):=\cF^{-1} f(x):=(2 \pi)^{-d/2}\int_{\mR^d} \e^{i\xi\cdot x}f(\xi)\dif\xi, \quad x\in\mR^d.
$$
For any $f\in\sS'(\mR^d)$,
\begin{align*}
\<\hat{f},\varphi\>:=\<f,\hat{\varphi}\>,\ \ \<\check{f},\varphi\>:=\<f,\check{\varphi}\>,\ \ \hbox{for } \forall\varphi\in\sS(\mR^d).
\end{align*}
Let $\chi:\mR^{d}\to[0,1]$ be a smooth radial function with
\begin{align*}
\chi(\xi)=
\begin{cases}
1, & \ \  |\xi|\leq 1,\\
0, &\ \ |\xi|>3/2.
\end{cases}
\end{align*}
Define $\psi(\xi):=\chi(\xi)-\chi(2\xi)$ and for $j\in\mN_0$,
\begin{align}\label{S2:Bscal}
\psi_j(\xi){:=}\psi(2^{-j}\xi).
\end{align}
Let $B_r := \{\xi\in \mR^d \mid  |\xi|\leq r\}$ for $r>0$. It is easy to see that $\psi\geq 0$,  supp$\psi\subset B_{3/2}/B_{1/2}$, and
\begin{align}\label{eq:PL01}
\chi(2\xi)+\sum_{j=0}^{k}\psi_j(\xi)=\chi(2^{-k}\xi)\to 1,\ \ \hbox{as}\ \ k\to\infty.
\end{align}
The block operators $\cR_j$ are defined on $\sS'(\mR^d)$ by
\begin{align*}
\cR_jf:=
\begin{cases}
\cF^{-1}(\chi\cF{f})=\check\chi* f, &j=-1,\\
\cF^{-1}(\psi_j\cF{f})=\check\psi_j* f, & j\geq 0.
\end{cases}
\end{align*}

\br
For $j\ge-1$,  by definitions, one sees that 
\begin{align}\label{S2:Rk25}
\cR_j=\cR_j\widetilde\cR_j,\quad \text{where $\widetilde\cR_j{:=}\sum_{\ell=-1}^1\cR_{j+\ell}$ with $\cR_{-2}:=0$},
\end{align}
and $\cR_j$ is symmetric in the sense of 
\begin{align}\label{S2:Rk26}
\int_{\mR^d}\cR_jf(x)g(x)\dif x=\int_{\mR^d}f(x)\cR_jg(x)\dif x, \ \ f \in \sS'(\mR^d),\, g \in \sS(\mR^d).
\end{align}
\er

Here is the definition of Besov spaces.
\bd[Besov spaces]\label{iBesov}
For any $\beta\in\mR$ and $p,q\in[1,\infty]$, the Besov space $\bB_{p,q}^\beta(\mR^d)$ is defined by
$$
\bB_{p,q}^\beta(\mR^d):=\Big\{f\in\sS'(\mR^d) \mid \|f\|_{\bB^\beta_{p,q}}:= \[ \sum_{j \geq -1}\left( 2^{\beta j} \|\cR_j f\|_{p} \right)^q \]^{1/q}  <\infty\Big\}.
$$
If $p=q=\infty$, it is in the sense
$$
\bB_{\infty,\infty}^\beta(\mR^d):=\Big\{f\in\sS'(\mR^d) \mid \|f\|_{\bB^\beta_{\infty,\infty}}:= \sup_{j \geq -1} 2^{\beta j} \|\cR_j f\|_{\infty} <\infty\Big\}.
$$
\ed

Recall the following Bernstein's inequality (cf. \cite{BCD11}*{Lemma 2.1}).

\bl[Bernstein's inequality]
For any $k\in\mN$, there is a constant $c=c(d,k)>0$ such that for all $j\ge-1$,
\begin{align*}
\|\nabla^k\cR_j f\|_\infty\le c 2^{kj}\|\cR_jf\|_\infty.
\end{align*}
In particular, for any $\alpha \in \mR$,
\begin{align}\label{S2:Bern}
\|\nabla^k f\|_{\bB^{\alpha}_{\infty,\infty}}\le c \|f\|_{\bB^{\alpha+k}_{\infty,\infty}}.
\end{align}
\el

\br[Equivalence between Besov spaces and H\"older  spaces]
If $\beta>0 \text{ and $\beta\notin\mN$}$, we have the  following equivalence between $\bB_{\infty,\infty}^\beta (\mR^d)$ and $ \bC^\beta  (\mR^d)$: (cf. \cite{Tr92})
\begin{align}\label{LJ1}
\|f\|_{\bB_{\infty,\infty}^\beta}\asymp\|f\|_{\bC^\beta}.
\end{align}
However, for any $n\in\mN_0$, we only have one side control that is
\begin{align}\label{eq:LJ1}
\|f\|_{\bB^n_{\infty,\infty}}\lesssim\|f\|_{\bC^n}.
\end{align}
By Bernstein's inequality, we have that for any $|h|<1/2$,
\begin{align*}
|f(x+h)-f(x)|
&\stackrel{\eqref{eq:PL01}}{\le} \sum_{j \ge-1}|\cR_j f(x+h)-\cR_jf(x)|
\lesssim \sum_{ j <-\log_2 |h|}\|f\|_{ \bB^1_{\infty,\infty}}|h|+\sum_{j\ge -\log_2 |h|}2^{-j}\|f\|_{ \bB^1_{\infty,\infty}}\\
& \lesssim \|f\|_{\bB^1_{\infty,\infty}}|h|(\log_2|h|^{-1}+1),
\end{align*}
and for any $|h|\geq 1/2$,
\begin{align*}
|f(x+h)-f(x)|&\le 2\|f\|_{\infty}\le 4 |h|\|f\|_{\bB^1_{\infty,\infty}}.
\end{align*}
Thus, by \eqref{S2:Bern}, we obtain that
\begin{align}\label{Pre:BH}
\sup_{x\ne y}\frac{|\nabla^k f(x)-\nabla^k f(y)|}{|x-y|(\log_2^+|x-y|^{-1}+1)}\lesssim\|f\|_{ \bB^{k+1}_{\infty,\infty}},\quad \text{for any $k\in\mN_0$}.
\end{align}
\er

Now we introduce the estimate of Littlewood-Paley's type for the heat kernel $p_\alpha(t,x)$. The same result is proved  in \cite{HWZ20}*{Lemma 3.1} for $\alpha=2$ and \cite{CHZ20}*{Lemma 3.3} and \cite{HWW20}*{Lemma 2.12} for $\alpha\in(0,2)$. For reader's convenience, we give a proof here.
\bl
Let $\alpha\in(0,2)$. There is a constant $c=c(\alpha,d)>1$ such that for all $j\ge-1$ and $T>0$,
\begin{align}\label{S2:EHK}
\int_0^T\int_{\mR^d}|\cR_j p_\alpha(t,x)|\dif x\dif t\le c(1+T)2^{-\alpha j}.
\end{align}
\el
\bpf
First of all, by the scaling property \eqref{eq:HG01}, we have that for any $m\in\mN_0$,
\begin{align}\label{eq:PL02}
\int_{\mR^d}|(\Delta^mp_\alpha)(t,x)|\dif x = t^{-2m/\alpha}\int_{\mR^d}|\Delta^m p_\alpha(1,x)|\dif x\lesssim  t^{-2m/\alpha}.
\end{align}
For $j=-1$, we have
\begin{align*}%\label{S2:EHK-1}
\int_0^T\int_{\mR^d}|\cR_{-1} p_\alpha(t,x)|\dif x\dif t\lesssim\int_0^T \|p_\alpha(t,\cdot)\|_{1}\dif t= T 2^{-\alpha }2^{\alpha }\leq T 2^{\alpha } .
\end{align*}
For $j \ge 0$,
by \eqref{S2:Bscal} and the change of variables,
\begin{align}\label{AZ01}
\int_{\mR^d}|\cR_{j} p_\alpha(t,x)|\dif x=2^{-jd}\int_{\mR^d}\Big|\int_{\mR^d}p_\alpha(t,2^{-j}(x-y))\check\psi(y)\dif y\Big|\dif x.
\end{align}
Notice that the support of $\psi$ is contained in an annulus. By \cite[(1.2.1)]{Gr14}, we have that $\Delta^{-m}\check \psi$ is a well-defined Schwartz function where
\begin{align*}
\cF(\Delta^{-m} \check{\psi})(\xi) :=(-|\xi|^2)^{-m} \psi(\xi)\in \sS(\mR^d),\, m\in \mN_0.
\end{align*}
Based on this, we have $\check \psi=\Delta^m\Delta^{-m}\check\psi,m\in \mN_0 $ and
\begin{align*}
\int_{\mR^d}p_\alpha(t,2^{-j}(x-y))\check\psi(y)\dif y
& =\int_{\mR^d}\Delta^m p_\alpha(t,2^{-j}(x-y))(\Delta^{-m}\check\psi)(y)\dif y.
\end{align*}
Hence, 
\begin{align*}
\int_{\mR^d}|\cR_{j} p_\alpha(t,x)|\dif x &  \stackrel{\eqref{AZ01}}{\lesssim} 2^{-jd} \int_{\mR^d}|\Delta^mp_\alpha(t,2^{-j} x)|\dif x\\
& = 2^{-2jm}\int_{\mR^d}|(\Delta^mp_\alpha)(t,x)|\dif x
 \stackrel{\eqref{eq:PL02}}{\lesssim} 2^{-2jm}t^{-2m/\alpha}.
\end{align*}
Then, considering the cases $m=0 $ and $m=2$, one sees that
\begin{align*}
\int_0^T\int_{\mR^d}|\cR_{j} p_\alpha(t,x)|\dif x\dif t&=\(\int_0^{2^{-\alpha j}}+\int_{2^{-\alpha j}}^T\)\int_{\mR^d}|\cR_{j} p_\alpha(t,x)|\dif x\dif t\\
&\lesssim \int_0^{2^{-\alpha j}}\dif t+2^{-4j}\int_{2^{-\alpha j}}^T t^{-4/\alpha}\dif t \lesssim 2^{-\alpha j}.
\end{align*}
The proof is finished.
\epf

Following the method used in \cites{HWZ20,HWW20}, we give a well-known   priori estimate of Besov-type by \eqref{S2:EHK}.  The result is seen as Schauder's estimate when $p=q=\infty$ in the literature. In the sequel, for a Banach space $\mB$ and $T>0$, $q\in[1,\infty]$, we denote by
$$
\mL_T^q(\mB):= L^q([0,T];\mB),\ \  \mL^q_T:=L^q([0,T]\times \mR^d).
$$

\bl\label{lem:PK01}
Let  $\alpha\in(0,2)$, $\beta\in\mR$. For any $p\in[1,\infty]$ and $q\in[1,\infty]$, there is a constant $c=c(d,\alpha,\beta,p,q)>0$ such that for all $(u,f)\in\sS'\times\sS'$ with
\begin{align*}
\p_t u=\Delta^{\alpha/2}u+f,\quad u(0)=u_0,
\end{align*}
in the following weak sense
\begin{align*}
\<u(t),\varphi\>=\<u_0,\varphi\>+\int_0^t\<u(s),\Delta^{\alpha/2}\varphi\>\dif s+\int_0^t\<f(s),\varphi\>\dif s,\ \ \forall\varphi\in\sS(\mR^d),
\end{align*}
and for any $T>0$,
\begin{align}\label{S2:Sch}
\|u\|_{ \mL^q_T(\bB^{\alpha+\beta}_{p,q})}\le c\(T^{1/q}\|u_0\|_{ \bB^{\alpha+\beta}_{p,q}}+(1+T)\|f\|_{ \mL^q_T(\bB^{\beta}_{p,q})}\).
\end{align}
\el

\bpf%[Proof of Lemma \ref{lem:PK01}]
We only give the proof under $q\in[1,\infty)$, since the case of $q=\infty$ is similiar and easier. Let $\{\rho_\eps\}_{\eps>0}$ be a usual modifier on $\mR^{d}$. Then $u_\eps:=u*\rho_\eps$ and $f_\eps=f*\rho_\eps$ satisfy
\begin{align*}
\p_t u_\eps(t,x)=\Delta^{\alpha/2}u_\eps(t,x)+f_\eps(t,x),\quad u_\eps(0)=u_0*\rho_\eps.
\end{align*}
Thus, without loss of generality, we assume that $u,f,u_0\in C^\infty$. For any $t\in[0,T]$, let $u^t(s):=u(t-s)$ and $f^t(s):=f(t-s)$ for any $s\in(0,t)$. Obviously,
\begin{align*}
\p_s u^t(s,x)+\Delta^{\alpha/2}u^t(s,x)=-f^t(s,x),\quad u^t(0)=u(t).
\end{align*}
By It\^o's formula (cf. \cite[Theorem 5.1]{IW89}), we have
\begin{align*}
\mE u^t(t,x+ L_t)=u^t(0,x)-\int_0^t\mE f^t(s,x+  L_s)\dif s.
\end{align*}
Then, we have Duhamel's formula:
\begin{align*}
u(t,x)=\int_{\mR^d}p_\alpha(t,x-y)u_0(y)\dif y+\int_0^t\int_{\mR^d}p_\alpha(s,x-y)f(t-s,y)\dif y\dif s.
\end{align*}
Taking $\cR_j$ for both sides, by \eqref{S2:Rk25} and \eqref{S2:Rk26}, we get
\begin{align*}
\cR_ju(t,x)=\int_{\mR^d}p_\alpha(t,x-y)\cR_ju_0(y)\dif y+\int_0^t\int_{\mR^d}\cR_jp_\alpha(s,x-y)\widetilde\cR_jf(t-s,y)\dif y\dif s.
\end{align*}
From this, by Minkowski's inequality and H\"older's inequality, one sees that
\begin{align*}
\|\cR_ju\|_{\mL^q_T(L^p)}&\le T^{1/q}\|\cR_ju_0\|_p+ \int_0^T\|\cR_jp_\alpha(s)\|_1\|\widetilde\cR_jf(\cdot-s) 1_{\cdot>s}\|_{\mL^q_T(L^p)}\dif s\\
&\stackrel{\eqref{S2:EHK}}{\lesssim} T^{1/q}\|\cR_ju_0\|_p+ (1+T)2^{-\alpha j}\|\widetilde\cR_jf \|_{\mL^q_T(L^p)}.
\end{align*}
By definitions and Fubini's theorem, we have
$$
\|u\|_{\mL_T^q(\bB^{\alpha+\beta}_{p,q})}^q
=\sum_{j\ge-1}2^{(\alpha+\beta)qj}\|\cR_ju\|_{\mL^q_T(L^p)}^q.
$$
Therefore,
\begin{align*}
\|u\|_{\mL_T^q(\bB^{\alpha+\beta}_{p,q})}^q
&\lesssim T \sum_{j\ge-1}2^{(\alpha+\beta)qj} \|\cR_ju_0\|_p^q+(1+T)^q\sum_{j\ge-1}2^{\beta qj}\|\cR_jf \|_{\mL^q_T(L^p)}^q\\
&\lesssim T\|u_0\|_{ \bB^{\alpha+\beta}_{p,q}}^q+(1+T)^q\|f\|_{ \mL^q_T(\bB^{\beta}_{p,q})}^q
\end{align*}
which implies the desired estimate.
\epf

\br[cf. {\cite[Section 3]{HWZ20}}]
The above result is true for $\alpha = 2$. In this case, we should cosider Brownian motion and { the Laplacian} $\Delta$ by the same way . 
\er
\br
Here we compare Schauder estimates in H\"older spaces, Besov spaces and Sobolev spaces.
By \eqref{S2:Sch} for $p=q=\infty$ and \eqref{LJ1}, we obtain the classical Schauder's estimate for $\alpha=2$:
\begin{align*}
\|u\|_{\mL^\infty_T(\bC^{2+\beta})}\lesssim \|u_0\|_{\bC^{2+\beta}}+\|f\|_{ \mL^\infty_T(\bC^{\beta})}, \quad \beta\in(0,1).
\end{align*}
It is well-known that Schauder's estimate is not true for $\beta=0$. But the lemma above tells us that
 \begin{align*}
\|u\|_{ \mL^\infty_T(\bB^{2}_{\infty,\infty})}\lesssim \|u_0\|_{ \bB^{2}_{\infty,\infty}}+\|f\|_{ \mL^\infty_T(\bB^{0}_{\infty,\infty})}.
\end{align*}
Furthermore, by \eqref{Pre:BH} and \eqref{eq:LJ1}, we get 
\begin{align*}
|\nabla u(t,x)-\nabla u(t,y)|\lesssim |x-y|\(1+\log_2^+|x-y|^{-1}\)\(\|u_0\|_{\bC^{2}}+\|f\|_{\mL^\infty_T}\).
\end{align*}
In Sobolev spaces,  it holds that
\begin{align*}
\|u\|_{\mL_T^q (W^{2,p})}\lesssim \|u_0\|_{\mL_T^q (W^{2,p})}+\|f\|_{\mL_T^q (L^{p})}
\end{align*}
with $p,q\ne1,\infty$ (see \cites{Kr01,XXZZ20} and  references therein). However, { $p,q=1$ or $\infty$ are} allowed in Besov case.
\er

\section{Estimates of heat kernels for Euler-Maruyama scheme}\label{Sec3}

In this section, assume that $\alpha \in (1,2)$ and $b : \mR_+\times \mR^d \to \mR^d$ is a bounded measurable function. Fix $T>0$ and $x\in \mR^d$. Consider the following Euler scheme $X_t^N(x)$: $X_0^N = x$, and 
\begin{align}\label{eq:WE01}
X_t^N = x  + \int_0^t b(s,X_{\phi_N(s)}^N) \dif s + L_t, \ \ t \in (0,T],
\end{align}
where $N \in \mN$, $\phi_N(s) := k h$ for $s\in [kh,(k+1)h]$ with $h := T/N$ and $k = 0,1,\cdots,N-1$. First of all, we prove the following Duhamel's formula for the Euler scheme. 

\bl[Duhamel's formula]\label{lem:AA01}
Let $\alpha\in(1,2)$. For each $t \in (0,T]$ and $x\in \mR^d$, $ X_t^N(x)$ admits a density $p_x^N(t,\cdot)$  satisfing the following Duhamel's formula:
\begin{align}\label{eq:BB02}
p_x^N(t,y) = p_\alpha(t,x-y) + \int_0^t \mE \[b(s, X_{\phi_N(s)}^N)\cdot \nabla p_\alpha(t-s,X_s^N-y)\]\dif s.
\end{align}
\el

\begin{proof}
Fix $t\in(0,T]$ and $f \in C_c^\infty(\mR^d)$. Letting $s \in [0,t]$ and 
$$
u(s,x): = p_\alpha(t-s,\cdot)*f(x)= \int_{\mR^d}p_\alpha(t-s,x-y) f(y) \dif y,
$$
by \eqref{eq:AA11}, it is easy to see that $u(s,x)$ solves the following equation:
\begin{align}\label{eq:BB01}
( \p_s + \Delta^{\alpha/2} ) u =0,\ \ u(t,x) = f(x).
\end{align}
By It\^o's formula (cf. \cite[Theorem 5.1]{IW89}), we have \begin{align*}
u(t,X_t^N) = \ \ 
& u(0,x) 
+ \int_0^t (\p_s u) (s,X_{s-}^N) \dif s 
+ \int_0^t b(s, X_{\phi_N(s)}^N) \cdot  \nabla  u (s,X_{s-}^N) \dif s\\
& + \int_0^{t} \int_{|z|>1} \( u (s,X_{s-}^{N}+z) - u (s,X_{s-}^{N}) \) N(\dif s, \dif z)\\
& +  \int_0^{t} \int_{0<|z|\leq 1} \( u (s,X_{s-}^{N}+z) - u (s,X_{s-}^{N}) \) \tilde N(\dif s, \dif z)\\
&+ \int_0^t \int_{\mR^d\setminus \{0\}} \( u (s,X_{s}^{N}+z\1_{|z|\leq 1}) - u (s,X_{s}^{N}) - z \1_{|z|\leq 1} \cdot \nabla u(s,X_{s}^{N}) \) \nu^{(\alpha)}( \dif z) \dif s.
\end{align*}
Observe that a c$\rm\grave{a}$dl$\rm\grave{a}$g function can have at most a countable number of  jumps. Taking the expectation for both sides in the above equality, by \cite[Section 3]{IW89},  \eqref{eq:AA05} and \eqref{eq:BB01}, we obtain that for any $f \in C_c^\infty(\mR^d)$,
\begin{align*}
\mE f(X_t^N) = \mE u(t,X_t^{N}) =  
u(0,x)
+  \int_0^t  \mE\( b(s, X_{\phi_N(s)}^N) \cdot \nabla u (s,X_{s}^N) \)\dif s.
\end{align*}
Furthermore, since
$$
\int_0^t \int_{\mR^d} |\nabla p_\alpha(s,y)| \dif y \dif s \stackrel{\eqref{eq:HG01}}{=} \int_{\mR^d}|\nabla p_\alpha(1,x)|\dif x\int_0^t s^{-1/\alpha} \dif s <\infty, \ \ \hbox{if } \alpha \in (1,2),
$$
we derive the desired Duhamel's formula.
\end{proof}

\br
For any { general} initial value $X_0^N = X_0 \in \sF_0$, since $L$ is independent of $X_0$, $X_t^N(x)$ defined by \eqref{eq:WE01} is independent of $X_0$. Consequently, by
\cite[Lemma 3.11]{Ka02}, the Euler scheme $ X_t^N$ with initial value $X_0$ also has a density $p_{X_0}^N(t,y)$ given by
\begin{align}\label{eq:EE01}
p_{X_0}^N(t,y) = \int_{\mR^d} p_x^N(t,y) \mP\circ X_0^{-1} (\dif x).
\end{align}
\er

The following uniform estimate for $p_x^N(t,y)$ was proved by Huang, Suo and Yuan \cite{HSY21} when the coefficient $b$ takes the form $b(x)$. For the convenience of readers, we show it again in the way of \cite{HRZ20}.

\bt\label{thm:AA01}
Let $\alpha\in(1,2)$. For any $T>0$, there is a constant $c = c(d,\alpha,T,\|b\|_\infty)>0$ such that for any $N \in \mN$, $t \in (0,T]$ and $x,y \in \mR^d$,
\begin{align}\label{S2:Ep}
p^N_x(t,y) \leq c p_\alpha(t,x-y).
\end{align}
\et

\begin{proof}%[Proof of Theorem \ref{thm:AA01}]
For the simplicity, we use a little confused notation $\|b\|_\infty:=\|b\|_{\mL^\infty_T}$ in the following. First of all, by \eqref{eq:AA09},  \eqref{eq:AA08} and \eqref{eq:AA02}, we know that there is a constant $ c_0=c_0(d,\alpha)>2$ such that
\begin{align}\label{eq:DD01}
|\nabla p_\alpha(t,x)| \leq c_0 t^{-1/\alpha} p_\alpha(t,x).
\end{align}
and 
\begin{align}\label{eq:DD02}
p_\alpha(t,x+z)\leq c_0 p_\alpha(t,x),\ \ \hbox{if}\ \  |z| \leq 2 t^{1/\alpha}.
\end{align}
Below, we fix this constant $c_0$ and $T>0$. Let $\eps>0$ be small enough such that
$$
\ell_\eps:= c_0^2\, \tfrac{\alpha}{\alpha-1} \|b\|_\infty\eps^{(\alpha-1)/\alpha}  \leq 1/2.
$$
Without loss of generality, we assume
$$
N \geq (T(\tfrac{1}{2}\|b\|_\infty)^{\alpha/(\alpha-1)} ) \vee (2T/\eps).
$$ 
Denote
\begin{align*}
h:=T/N \ \ \text{and} \ \ M:=[\eps/h]\in\mN.
\end{align*}
 Then, we have $\|b\|_\infty \leq 2 h^{-1+1/\alpha}$ and $\eps>h$.

\medskip\noindent
{\bf (Step 1)} In this step, by induction, we prove the following result: for $k = 1,2,\ldots,M\wedge N$, 
\begin{align}\label{eq:DD03}
p_x^N(kh,y) \leq c_0 p_\alpha(kh,x-y).
\end{align}
For $k=1$, noting that $ X_h^N = x + \int_0^h b(s,x)\dif s + L_h$ with $\|b\|_\infty \leq 2 h^{-(\alpha-1)/\alpha}$, by \eqref{eq:DD02} we get that
$$
p_x^N(h,y) = p_\alpha(h,y - x -  \int_0^h b(s,x)\dif s ) \leq  c_0 p_\alpha(h,x-y).
$$
Suppose now that \eqref{eq:DD03} holds for $j = 1,2,\ldots,k-1$. By Duhamel's formula \eqref{eq:BB02}, we see that
\begin{align}\label{eq:DD04}
p_x^N(kh,y) - p_\alpha(kh,x-y) 
&= \int_0^{kh} \mE\[ b(s,X_{\phi_{N}(s)}^N)\cdot \nabla p_\alpha(kh-s,X_s^N-y) \]\dif s\nonumber\\
& = \sum_{j=0}^{k-1} \int_{jh}^{(j+1)h} I_j^N(s)  \dif s,
\end{align}
where
$
I_j(s) :=  \mE\[ b(s,X_{jh}^N)\cdot \nabla p_\alpha (kh-s,X_s^N-y)\].
$
Observe that for $s\in (jh,(j+1)h)$,
$$
X_{s}^N = X_{jh}^N + \int_{jh}^sb(r,X_{jh}^N)\dif r+ (L_s - L_{jh}).
$$
Since $L_s  - L_{jh}$ is independent of $X_{jh}^N$ and has density $p_\alpha(s-jh,\cdot)$, 
by \cite[Lemma 3.11]{Ka02} and C-K equations \eqref{eq:CC00}, we have 
\begin{align*} 
I_j(s) 
& = \mE \[ b(s,X_{jh}^N)\cdot \nabla p_\alpha (kh-s)*p_\alpha(s-jh)\(X_{jh}^N + \int_{jh}^sb(r,X_{jh}^N)\dif r -y\) \]\\
& = \mE \[ b(s,X_{jh}^N)\cdot \nabla p_\alpha \(kh-jh, X_{jh}^N + \int_{jh}^sb(r,X_{jh}^N)\dif r -y\) \]\\
& \leq \|b\|_\infty \int_{\mR^d}  | \nabla p_\alpha | \(kh-jh, z -y+ \int_{jh}^sb(r,z)\dif r \) p_x^N(jh,z) \dif z.
\end{align*}
Furthermore, by \eqref{eq:DD01}, \eqref{eq:DD02} and induction hypothesis, we obtain that for $s\in(jh,(j+1)h)$,
\begin{align*}
I_j(s) 
& \leq \|b\|_\infty  (kh-jh)^{-1/\alpha} c_0^2 \int_{\mR^d} p_\alpha(kh-jh,z-y) \cdot c_0 p_\alpha(jh,x-z) \dif z\\
& \leq c_0\, \tfrac{\alpha-1}{\alpha} \ell_\eps \eps^{-(\alpha-1)/\alpha}(kh-s)^{-1/\alpha}  p_\alpha(kh,x-y),
\end{align*}
where we have used $h\|b\|_\infty \leq 2 h^{1/\alpha}$. Substituting this into \eqref{eq:DD04}, we get, since $kh \leq Mh\leq \eps$ and $\alpha\in(1,2)$, that
\begin{align*}
|p_x^N(kh,y) - p_\alpha(kh,x-y) |
& \leq c_0\, \ell_\eps \eps^{-(\alpha-1)/\alpha}  p_\alpha(kh,x-y)  \tfrac{\alpha-1}{\alpha}\int_0^{kh} (kh-s)^{-1/\alpha} \dif s\\
& = c_0\, \ell_\eps \eps^{-(\alpha-1)/\alpha} (kh)^{(\alpha-1)/\alpha} p_\alpha(kh,x-y)\\
& \leq c_0\, \ell_\eps p_\alpha(kh,x-y),
\end{align*}
which implies that 
$$
p_x^N(kh,y) \leq (  c_0\, \ell_\eps+ 1 )  p_\alpha(kh,x-y) \leq  c_0\,p_\alpha(kh,x-y).
$$

\medskip\noindent
{\bf (Step 2)} Next we assume that $M < N$. Since $\phi_N(s+Mh) = \phi_N(s)+Mh$, we have
\begin{align*}
X^N_{t+Mh} 
& = X_{Mh}^N +\int_{Mh}^{t+Mh} b (s, X_{\phi_N(s)}^N) \dif s + (L_{t+Mh}  - L_{Mh} )\\
& = X_{Mh}^N +\int_{0}^{t} b (s+Mh, X_{\phi_N(s)+Mh}^N) \dif s + (L_{t+Mh}  - L_{Mh} ).
\end{align*}
For $t \in [0,Mh]$, letting
$$
\tilde X^N_{t}  =  X^N_{t+Mh},\ \  \tilde L_t  = L_{t+Mh}  - L_{Mh} ,
$$
we have
$$
\tilde X^N_{t} = \tilde X^N_{0} + \int_{0}^{t} b (s+ Mh, \tilde X_{\phi_N(s)}^N) \dif s +  \tilde L_t.
$$
Noting that $(\tilde L_t )_{t\geq 0} \stackrel{d}{=} (L_ t)_{t\geq 0}$, denoting by $\tilde p_z^{N}(t,\cdot)$ the density of $\tilde X^N_{t}$ with $\tilde X^N_{0} =z$, by  Step 1, we have
$$
\tilde p_z^{N}(jh,y) \leq c_0\, p_\alpha(jh,z-y),\ \ j = 1,\ldots,M.
$$
Hence, for $j = 1,\ldots,M$, by \eqref{eq:EE01}, \eqref{eq:DD03} and C-K equations \eqref{eq:CC00}, we obtain that
\begin{align*}
p_x^N((j+M)h,y) 
& = \int_{\mR^d} \tilde p_z^N(jh,y) p_x^N(Mh,z)\dif z\\
& \leq c_0^2 \int_{\mR^d} p_\alpha(jh,z-y) p_\alpha(Mh,x-z) \dif z\\
& = c_0^2 p_\alpha ((j+M)h,x-y),
\end{align*}
that is 
$$
p_x^N(kh,y) \leq c_0^2 p_\alpha (kh,x-y),\ \ k = M+1,\ldots,2M.
$$
Repeating the above procedure $[N/M]$-times, we get that 
$$
p_x^N(kh,y) \leq c_0^{[2T/\eps]+1} p_\alpha (kh,x-y),\ \ k = 1,\ldots,N.
$$
We point that the constant $c_0^{[2T/\eps]+1}$ is independent of $N$.

\medskip\noindent
{\bf (Step 3)} Observe that for  $t \in (kh,(k+1)h)$,
$$
X_t^N = X_{kh}^N + \int_{kh}^tb (s, X_{kh}^N)\dif s + (L_{t} - L_{kh}),
$$
where $L_{t} - L_{kh}$ is independent of $X_{kh}^N$. Thus, by \cite[Lemma 3.11]{Ka02} and \eqref{eq:DD02},
\begin{align*}
p_x^N(t,y) 
& = \int_{\mR^d} p_x^N(kh,z) p_\alpha(t-kh, z+\int_{kh}^tb(s,z)\dif s-y) \dif z\\
& \leq c_0^{[2T/\eps]+2} \int_{\mR^d} p_\alpha(kh,x-z) p_\alpha(t-kh, z-y) \dif z\\
&= c_0^{[2T/\eps]+2} p_\alpha(t, x-y).
\end{align*}
Here, we have used $h\|b\|_\infty \leq 2 h^{1/\alpha}$ and C-K equations \eqref{eq:CC00}.
\end{proof}

The following corollary is a combination of Theorem \ref{thm:AA01}, Lemma \ref{lem:AA02} and Lemma \ref{lem:AA03}.

\bc\label{cor:AA01}
Let $\mu_0(\dif x) = \mP\circ X_0^{-1}(\dif x)$ be the distribution of $X_0$ and $\alpha \in (1,2)$.
\begin{enumerate}[\quad(i)]
\item For any $T>0$, there is a constant $c= c(d,\alpha,T,\|b\|_\infty)>0$ such that for all $N \in \mN$, $t \in (0,T]$ and $x\in \mR^d$,
\begin{align}\label{Se3:Main}
p_{X_0}^N(t,y) \leq c \int_{\mR^d} p_\alpha(t,x-y) \mu_0(\dif x).
\end{align}
\item For any $T>0$ and $\beta \in (0,\alpha-1)$, there is a constant $c=c(d,\alpha,T,\|b\|_\infty,\beta)>0$ such that for all $N\in \mN$, $t \in(0,T]$ and $y_1,y_2 \in \mR^d$,
\begin{align}\label{S2:Holder}
|p_{X_0}^N(t,y_2) -  p_{X_0}^N(t,y_1) |  \leq c |y_2 - y_1|^{\beta}  t^{-\beta/\alpha} \sum_{i=1,2} \int_{\mR^d} p_\alpha(t,x-y_i) \mu_0(\dif x).
\end{align}
\item For any $T>0$ and $\beta \in (0,\alpha-1)$, there is a constant $c=c(d,\alpha,T,\|b\|_\infty,\beta)>0$ such that for all $N\in \mN$, $t_1,t_2\in(0,T]$ and $y \in \mR^d$,
$$
|p_{X_0}^N(t_2,y) -  p_{X_0}^N(t_1,y) |  \leq c |t_2 - t_1|^{\beta/\alpha} \sum_{i=1,2} t_i^{-\beta/\alpha} \int_{\mR^d} p_\alpha(t_i,x-y) \mu_0(\dif x).
$$
\end{enumerate}
\ec

\begin{proof}
$(i)$ is a direct consequence of \eqref{eq:EE01} and Theorem \ref{thm:AA01}. 

$(ii)$
By Duhamel's formula \eqref{eq:BB02} and \eqref{eq:EE01}, we have
\begin{align*}
|p_{X_0}^N(t,y_2) - &  p_{X_0}^N(t,y_1) | 
\leq  \sI_1 +   \sI_2,
\end{align*}
where
$$
\sI_1 := \int_{\mR^d}| p_\alpha(t,x-y_2) -  p_\alpha(t,x-y_1) | \mu_0(\dif x),
$$
and
$$
 \sI_2 := \|b\|_\infty  \int_{0}^{t}  \int_{\mR^d} |\nabla p_\alpha(t-s,y_1-z) - \nabla p_\alpha(t-s,y_2-z) | p_{X_0}^N(s,z) \dif z\dif s,
$$
For $\sI_1$, by \eqref{eq:AA01}, we have
$$
\sI_1 \lesssim |y_2 - y_1|^{\beta}  t^{- \beta/\alpha} \sum_{i=1,2} \int_{\mR^d} p_\alpha(t,x-y_i) \mu_0(\dif x).
$$
For $\sI_2$, by \eqref{eq:AA01}, $(i)$ and C-K equations \eqref{eq:CC00}, we obtain that
\begin{align*}
\sI_2 
& \lesssim |y_2 - y_1|^{\beta} 
\int_0^{t} (t-s)^{-(1+\beta)/\alpha} \sum_{i=1,2} 
\left( \int_{\mR^d} p_\alpha(t - s,z-y_i ) \[\int_{\mR^d} p_\alpha(s,x-z) \mu_0(\dif x)\] \dif z \right) 
\dif s\\
& = |y_2 - y_1|^{\beta}  
\int_0^{t} (t -s)^{-(1+\beta)/\alpha} \dif s 
\sum_{i=1,2} 
\int_{\mR^d} p_\alpha(t,x-y_i) \mu_0(\dif x)\\
& \lesssim |y_2 - y_1|^{\beta} t^{(\alpha-1-\beta)/\alpha}
\sum_{i=1,2}
\int_{\mR^d} p_\alpha(t,x-y_i) \mu_0(\dif x),
\end{align*}
where we have used $\beta \in (0,\alpha-1)$. 

$(iii)$ Suppose that $t_1<t_2$. By Duhamel's formula \eqref{eq:BB02} and \eqref{eq:EE01}, we have
\begin{align*}
|p_{X_0}^N(t_2,y) - &  p_{X_0}^N(t_1,y) | 
\leq  \sJ_1 +   \sJ_2 +  \sJ_3,
\end{align*}
where
$$
\sJ_1 := \int_{\mR^d}| p_\alpha(t_2,x-y) -  p_\alpha(t_1,x-y) | \mu_0(\dif x),
$$
$$
 \sJ_2 := \|b\|_\infty  \int_{t_1}^{t_2}  \int_{\mR^d} |\nabla p_\alpha(t_2-s,z-y)| p_{X_0}^N(s,z) \dif z,
$$
and
$$
 \sJ_3 :=\|b\|_\infty   \int_{0}^{t_1}  \int_{\mR^d} |\nabla p_\alpha(t_2-s,z-y) - \nabla p_\alpha(t_1-s,z-y) | p_{X_0}^N(s,z) \dif z
$$
For $\sJ_1$, by \eqref{eq:AA03}, we have
$$
\sJ_1 \lesssim |t_2 - t_1|^{\beta/\alpha} \sum_{i=1,2} t_i^{-\beta/\alpha} \int_{\mR^d} p_\alpha(t_i,x-y) \mu_0(\dif x).
$$
For $\sJ_2$, by \eqref{eq:AA08}, $(i)$ and C-K equations \eqref{eq:CC00}, we get
\begin{align*}
\sJ_2 
& \lesssim \int_{t_1}^{t_2}  (t_2-s)^{-1/\alpha} \left( \int_{\mR^d} p_\alpha(t_2-s,z-y) \int_{\mR^d} p_\alpha(s,x-z) \mu_0(\dif x) \dif z\right) \dif s\\
& = \int_{t_1}^{t_2}   (t_2-s)^{-1/\alpha} \dif s \int_{\mR^d}  p_\alpha(t_2,x-y) \mu_0(\dif x)\\
& \lesssim  (t_2 - t_1)^{-1/\alpha+1} \int_{\mR^d}  p_\alpha(t_2,x-y) \mu_0(\dif x).
\end{align*}
Since $\beta \in (0,\alpha-1)$, we have
$$
0\leq  (t_2 - t_1)^{-1/\alpha+1} \leq |t_2 - t_1|^{\beta/\alpha} t_2^{-\beta/\alpha+(\alpha-1)/\alpha} \leq |t_2 - t_1|^{\beta/\alpha}  t_2^{-\beta/\alpha} T^{(\alpha-1)/\alpha}.
$$
Hence,
$$
\sJ_2 \lesssim   |t_2 - t_1|^{\beta/\alpha}  t_2^{-\beta/\alpha} \int_{\mR^d}  p_\alpha(t_2,x-y) \mu_0(\dif x).
$$
For $\sJ_3$, by \eqref{eq:AA03}, $(i)$ and C-K equations \eqref{eq:CC00}, we obtain that
\begin{align*}
\sJ_3 
& \lesssim |t_2 - t_1|^{\beta/\alpha} \sum_{i=1,2} \int_0^{t_1} (t_i -s)^{-(1+\beta)/\alpha} \left( \int_{\mR^d} p_\alpha(t_i - s,z-y) \int_{\mR^d} p_\alpha(s,x-z) \mu_0(\dif x) \dif z \right) \dif s\\
& = |t_2 - t_1|^{\beta/\alpha} \sum_{i=1,2} \int_0^{t_1} (t_i -s)^{-(1+\beta)/\alpha} \dif s  \int_{\mR^d} p_\alpha(t_i,x-y) \mu_0(\dif x)\\
& \leq{\tfrac{\alpha }{\alpha-1-\beta}} T^{(\alpha-1)/\alpha}
|t_2 - t_1|^{\beta/\alpha} \sum_{i=1,2} t_i^{-\beta/\alpha} \int_{\mR^d} p_\alpha(t_i,x-y) \mu_0(\dif x),
\end{align*}
where we have used $\beta \in (0,\alpha-1)$ and $0\le t_1<t_2\leq T$.

Combining the above calculations, we get the desired estimate.
\end{proof}

\section{Proof of Theorem \ref{in:Main}}\label{Sec:Main}

Let $(\Omega,\sF,(\sF_t)_{t\geq 0},\mP)$ be a complete filtered probability space, $L_t$ a $d$-dimensional symmetric and rotationally invariant $\sF_t$-adapted $\alpha$-stable process with $\alpha \in (1,2)$, $X_0$ an $\sF_0$-measurable random variable with distribution $\mu_0$. Let $T>0$, $N \in \mN$ and $h:=T/N$. Let $X_t^N$ be the Euler approximation of DDSDE \eqref{eq:AB} constructed in the introduction. From the construction, it is easy to see that $X_t^N$ solves the following SDE:
\begin{align}\label{eq:AC01}
X_t^N = X_0 + \int_0^t b^N(s,X_{\phi_N(s)}^N) \dif s + L_t,
\end{align}
where 
\begin{align*}
b^N(s,x) = \1_{\{s\geq h\}} b(s,x,\rho_{\phi_N(s)}^N(x))
\end{align*}
and 
\begin{align*}
\phi_N(s) = \sum_{j=0}^\infty jh \1_{[jh,(j+1)h)}(s).
\end{align*}
Trivially, $s-h \leq \phi_N(s)  \leq s$.

\medskip 

Let  $\mD$ be the space of all c$\rm\grave{a}$dl$\rm\grave{a}$g functions from $[0,T]$ to $\mR^d$. In the following, $\mD$ is equipped with Skorokhod topology which makes $\mD$ into a Polish space, and use $d_\mD$ to denote the associated metirc.

\bl\label{lem:CC01}
The sequence of laws for $(X_{\cdot}^N) $ in $(\mD,d_\mD)$ is tight.
\el

\begin{proof}
It is trivial that the sequence of distributions for $(X_0^N,L_0) \equiv (X_0,0)$ is tight in $\mR^d\times \mR^d$.
Taking $q \in (\alpha/2,\alpha)$, by Chebyshev's inequality, \eqref{eq:AC01} and the fact (cf. \cite[Lemma 2.4]{CZ18})
\begin{align}\label{eq:LO01}
\mE |L_t - L_s|^q \lesssim  |t-s|^{q/\alpha},\ \ q\in(0,\alpha),
\end{align}
 we obtain that for any $N \in \mN$, $R>0$ and $0\leq s<r<t \le T$,
\begin{align}\label{S4:NX01}
\begin{split}
\quad &\mP\( |X_r^N-X_s^N|\geq R, |X_t^N-X_r^N|\geq R \) \\  
\leq  \,& \mP\(  |L_r-L_s| + (r-s)\|b\|_\infty \geq R \)
 \times \mP\(|L_r-L_s|+ (t-r)\|b\|_\infty\geq R \) \\
\lesssim&   (r-s)^{q/\alpha} (t-r)^{q/\alpha} R^{-2q}
\leq (t-s)^{2q/\alpha} R^{-2q}.
\end{split}
\end{align}
Similarly, we have
$$
\lim_{\delta \downarrow 0} \sup_{N} \mP \(\ |X_\delta^N - X_0^N | \geq \eps\)  = 0,\ \  \forall \eps>0.
$$
Hence, combining the above calculations, by \cite[Theorem 4.1, p.\,355]{
JS03}, we see that the sequence $(X_{\cdot}^N)$ is tight.
\end{proof}

Let $ p_x^N(t,\cdot)$ be the distributional density of the Euler scheme $X_t^N(x)$ of SDE \eqref{eq:AC01} starting from $x$ at time $0$. Since for each $x\in\mR^d$, $X_t^N(x)$ is independent of $X_0$, the distributional density $\rho_t^N(\cdot)$ of $X_t^N$ with initial distribution $\mu_0$ is given by
\begin{align}\label{eq:RR01}
\rho_t^N(y) = \int_{\mR^d} p_x^N(t,y) \mu_0(\dif x).
\end{align}
Furthermore, by Theorem \ref{thm:AA01}, we have that for $q>1$,
\begin{align}\label{eq:LI01}
\begin{split}
\( \int_{\mR^d}  | \rho_{\phi_N(s)}^N(y) |^q \dif y \)^{1/q}
& = \( \int_{\mR^d}  \left| \int_{\mR^d} p_x^N(\phi_N(t),y) \mu_0(\dif x) \right|^q \dif y \)^{1/q}\\
& \lesssim
 \( \int_{\mR^d}   \int_{\mR^d} \left| p_\alpha(\phi_N(t),x-y)\right|^q \mu_0(\dif x)  \dif y \)^{1/q}
\stackrel{\eqref{eq:HG01}}{\lesssim} \phi_N(t)^{-\frac{d}{\alpha p} },
\end{split}
\end{align}
where $1/q+1/p=1$.

\bl\label{lem:AM}
For fixed $T>0$, there are a subsequence $(N_k
)_{k \in \mN}$ and a continuous function $\rho \in C((0,T] \times \mR^d)$ such that for any $M \in \mN$ with $M>1/T$,
\begin{align}\label{eq:DD04-1}
\lim_{k\to \infty} \sup_{|y|\leq M} \sup_{1/M\leq t\leq T} |\rho_t^{N_k} (y) - \rho_t(y)| = 0.
\end{align}
\el

\begin{proof}
By Theorem \ref{thm:AA01} and \eqref{eq:AA09}, we have that
$$
\sup_{|y|\leq M} \sup_{1/M\leq t\leq T} |\rho_t^N(y)|
\leq c \int_{\mR^d} \sup_{|y|\leq M} \sup_{1/M\leq t\leq T} |p_\alpha(t,x-y)|  \mu_0(\dif x) 
\leq c_M,
$$
where $c_M$  is independent of $N$. Moreover, by Corollary \ref{cor:AA01}, we have for any $\beta \in (0,\alpha-1)$, $t_1,t_2 \in [1/M,T]$ and $y_1,y_2 \in \mR^d$,
\begin{align}\label{eq:DF}
|\rho_{t_1}^N(y_1) -  \rho_{t_2}^N(y_2) |
\leq \,& |\rho_{t_1}^N(y_1) -  \rho_{t_2}^N(y_1) | + |\rho_{t_2}^N(y_1) -  \rho_{t_2}^N(y_2)| \nonumber\\
\lesssim \, & |t_1-t_2|^{\beta/\alpha} M^{\beta/\alpha} 
\sum_{i=1,2}
\int_{\mR^d} p_\alpha(t_i, x-y_1) \mu_0(\dif x) \nonumber\\
& + |y_1-y_2|^{\beta} M^{\beta/\alpha}
\sum_{i=1,2}
\int_{\mR^d} p_\alpha(t_2, x-y_i) \mu_0(\dif x)\nonumber\\
\stackrel{\eqref{eq:AA09}}{\lesssim}\, & 
M^{(d+\beta)/\alpha}(|t_1-t_2|^{\beta/\alpha}+ |y_1-y_2|^{\beta}),
\end{align}
where the implicit constants in the above $\lesssim$ are independent of $N$. Thus, by Ascolli-Arzela's theorem, we conclude the proof and have \eqref{eq:DD04-1}.
\end{proof}

Now we are in a position to give

\begin{proof}[Proof of Theorem \ref{in:Main}] 
{\bf(Existence)}
Fix $T>0$. For the simiplicity, we use a little confused notation $\|\cdot\|_\infty:=\|\cdot\|_{\mL^\infty_T}$ in some places. Let $\mQ_N$ be the law of $(X^N,L)$  in $\mD \times \mD$. By Lemma \ref{lem:CC01}, $\mQ_N$ is tight. Therefore, by Prokhorov's theorem (cf.  \cite[Theorem 16.3]{Ka02}), for the subsequence in Lemma \ref{lem:AM}, there are a subsubsequence $(n_j)_{j\geq 1}$ and a probability measure $\mQ$ on $\mD\times \mD$ so that
$$
\mQ_{n_j} \to \mQ \ \ \hbox{weakly}.
$$
Below, for simplicity of notations, we still denote the above subsequence by $\mQ_N$, $N \in \mN$. Then, by Skorokhod's representation theorem (cf. \cite[Theorem 4.30]{Ka02}), there are a probability space $(\widetilde \Omega,\widetilde \sF, \widetilde \mP)$ and random variables $\widetilde X, \widetilde L$ thereon such that
\begin{align}\label{eq:AG}
(\widetilde X^N, \widetilde L^N) \stackrel{}{\rightarrow} (\widetilde X , \widetilde L ),\ \ \widetilde \mP-\hbox{a.s.}
\end{align}
and
\begin{align}\label{eq:AF}
\widetilde \mP \circ (\widetilde X^N, \widetilde L^N)^{-1} = \mQ_N =  \mP \circ (X^N, L)^{-1},\ \ \widetilde \mP \circ  (\widetilde X , \widetilde L )^{-1} = \mQ.
\end{align}
In particular, the distributional density of $\widetilde X_t^N$ is $\rho_t^N$. Moreover, by Lemma \ref{lem:AM} and \eqref{eq:AG}, for any $t \in (0,T)$ and $\varphi \in C_c^\infty(\mR^d)$,
$$
\mE \varphi(\widetilde X_t) = \lim_{N \to \infty} \mE \varphi(\widetilde X_t^N) = \lim_{N \to \infty} \int_{\mR^d} \varphi(z) \rho_t^N(z) \dif z = \int_{\mR^d} \varphi(z) \rho_t(z) \dif z.
$$
In other words, $\rho_t$ is the density of $\widetilde X_t$. Define $\widetilde \sF^N_t := \sigma\{\widetilde X_s^N , \widetilde L_s^N ;s\leq t \}$. Noting that 
$$
\mP [L_t - L_s \in \cdot \mid \sF_s] = \mP \{ L_t - L_s \in \cdot \},
$$
we have
$$
\widetilde \mP [\widetilde L^N_t - \widetilde L^N_s \in \cdot \mid \widetilde \sF^N_s] = \mP \{ \widetilde L^N_t - \widetilde L^N_s \in \cdot \},
$$
which means that $\widetilde L^N_t$ is an $(\widetilde \sF^N_s)$-adapted $\alpha$-stable L\'evy process. Thus, by \eqref{eq:AC01} and \eqref{eq:AF} we obtain
\begin{align}\label{S3:SDEN}
\widetilde X_t^N = \widetilde X_0^N + \int_0^t b^N(s, \widetilde X^N_{\phi_N(s)})\dif s + \widetilde L_t^N,
\end{align}
where $b^N(s, \widetilde X^N_{\phi_N(s)}) 
= \1_{\{s\geq h\}} b(s ,\widetilde X^N_{\phi_N(s)},\rho_{\phi_N(s)}^N(\widetilde X^N_{\phi_N(s)}))$. We claim that
\begin{align}\label{S3:CoB}
\int_0^t b^N(s, \widetilde X^N_{\phi_N(s)})\dif s \to \int_0^t b(s , \widetilde X_{s}, \rho_s(\widetilde X_{s}))\dif s,
\end{align}
in probability as $N \to \infty$. Recalling the results in \cite[p.\,339]{JS03} and \eqref{eq:AG} , one sees that for $\tilde{P}$-a.s. $\omega$,  if $\Delta \widetilde X_t(\omega) =\Delta \widetilde L_t(\omega) =0 $, then
$$
\widetilde X_t^N(\omega) \to \widetilde X_t(\omega).
$$
Then, through taking $N\to\infty$ in \eqref{S3:SDEN}, it holds that for $\widetilde \mP$-a.s $\omega$,
\begin{align*}
\widetilde X_t(\omega) = \widetilde X_0(\omega) + \int_0^t b(s, \widetilde X_s(\omega),\rho_s(\widetilde X_s(\omega)))\dif s + \widetilde L_t(\omega),\quad t\in D_\omega,
\end{align*}
where
\begin{align*}
D_\omega:=\{t\in\mR_+ \mid \Delta \widetilde X_t(\omega)=\Delta\widetilde L_t(\omega)=0 \}.
\end{align*}
 Since $\widetilde X$ and $\widetilde L$ belong to $\mD$, $D_\omega^c$ is a countable set in $\mR_+$ and 
  \begin{align*}
\widetilde X_t(\omega) = \widetilde X_0(\omega) + \int_0^t b(s, \widetilde X_s(\omega),\rho_s(\widetilde X_s(\omega)))\dif s + \widetilde L_t(\omega),\quad t\in \mR_+,
\end{align*}
which derives the existence.

%\medskip 

Let us now prove \eqref{S3:CoB}. Indeed, observe that
\begin{align*}
 \mE \left| \int_0^t b^N(s, \widetilde X^N_{\phi_N(s)})\dif s- \int_0^t b(s , \widetilde X_{s}, \rho_s(\widetilde X_{s}))\dif s \right|
\leq 
 \sJ_1^N + \sJ_2^N +  T  \|b\|_\infty /N,
\end{align*}
where
\begin{align*}
 \sJ_1^N := \mE \int_h^t \Big|
b(s ,\widetilde X^N_{\phi_N(s)},\rho_{\phi_N(s)}^N(\widetilde X^N_{\phi_N(s)})) - b(s ,\widetilde X^N_{\phi_N(s)},\rho_{s}(\widetilde X^N_{\phi_N(s)})) \Big| \dif s
\end{align*}
and
\begin{align*}
\sJ_2^N := \mE \int_h^t \Big |b(s ,\widetilde X^N_{\phi_N(s)},\rho_{s}(\widetilde X^N_{\phi_N(s)})) - b(s , \widetilde X_{s}, \rho_s(\widetilde X_{s})) \Big| \dif s.
\end{align*}
\begin{enumerate}[(1)]
\item For $\sJ_1^N$, we have
\begin{align*}
\sJ_1^N
\leq &\, \mE \int_h^t \1_{ \{ |\widetilde X_{\phi_N(s)}^N|\leq R \} }  \big| b(s ,\widetilde X^N_{\phi_N(s)},\rho_{\phi_N(s)}^N(\widetilde X^N_{\phi_N(s)})) - b(s ,\widetilde X^N_{\phi_N(s)},\rho_{s}(\widetilde X^N_{\phi_N(s)})) \Big|  \dif s\\
& + 2 \|b\|_\infty  \int_h^t \widetilde \mP \( |\widetilde X_{\phi_N(s)}^N|> R \) \dif s
:=  \sJ_{11}^N(R) + \sJ_{12}^N(R) .
\end{align*}
Since
$$
| \rho_{\phi_N(s)}^N (x) - \rho_s (x) | \leq | \rho_{\phi_N(s)}^N (x) - \rho_s^N (x) | + | \rho_{s}^N (x) - \rho_s (x) |,
$$
by \eqref{eq:DD04-1} and \eqref{eq:DF}, we see that for each fixed $(s,x) \in (0,T]\times\mR^d$,
\begin{align*}
\lim_{N \to \infty}\1_{\{s\geq h\}} | \rho_{\phi_N(s)}^N (x) - \rho_s  (x) | = 0,
\end{align*}
which implies that for any $(s,x)\in\mR_+\times\mR^d$, by \eqref{eq:PO} we have
\begin{align}\label{S3:BB}
\lim_{N\to\infty} | b(s,x,\rho_{\phi_N(s)}^N(x)) - b(s,x,\rho_s(x))| = 0.
\end{align}
Moreover, by H\"older's equality and \eqref{eq:LI01}, we get
\begin{align}\label{eq:ZK01}
\begin{split}
\sJ^N_{11}(R)=&\int_h^t \int_{|x|\le R} \Big|b(s,x,\rho_{\phi_N(s)}^N(x))-b(s,x,\rho_s(x))\Big|\rho_{\phi_N(s)}^N(x)\dif x\dif s\\
\leq &\, \left[\int_h^t \int_{|x|\le R} \Big|b(s,x,\rho_{\phi_N(s)}^N(x))-b(s,x,\rho_s(x))\Big|^p\dif x \dif s \right]^{1/p}\\
& \qquad \times 
\left[\int_h^t \int_{|y|\leq R}  | \rho_{\phi_N(s)}^N(y) |^q \dif y  
\dif s \right]^{1/q}\\
\lesssim &\, \left[\int_h^t \int_{|x|\leq R} \Big|b(s,x,\rho_{\phi_N(s)}^N(x))-b(s,x,\rho_s(x))\Big|^p\dif x \dif s \right]^{1/p}\\
& \qquad \times 
\left[ \int_h^T (s-h)^{-\frac{d}{\alpha }(q-1) } \dif s \right]^{1/q}
\end{split}
\end{align}
provided by $1<q<1+{\alpha}/{d}$ and  $1/p+1/q=1$. Note that the implicit constant in the above $\lesssim$ is independent of $N,R$. Thus, for any $R >0$, by the dominate convergence theorem and \eqref{S3:BB}, we get that
\begin{align}\label{eq:KJ01}
\lim_{N \to \infty} \sJ^N_{11}(R)  = 0.
\end{align}
For $\sJ^N_{12}(R)$,  by \eqref{eq:AC01} , \eqref{eq:LO01} and  Chebyshev's inequality, we have
\begin{align}\label{eq:KJ02}
 \int_0^t \widetilde \mP  \( |\widetilde X_{\phi_N(s)}^N|> R \) \dif s
& =  \int_0^t \mP \( |X_{\phi_N(s)}^N|> R \) \dif s
\nonumber\\
& \lesssim T  \mP \( |X_0| + T \|b\|_\infty >R/2\)  + \int_0^t \frac{(\phi_N(s))^{1/\alpha}}{R/2} \dif s\nonumber\\
& \leq T  \mP \( |X_0| + T \|b\|_\infty >R/2\)  + T^{(\alpha+1)/\alpha} (R/2)^{-1}
\end{align}
which converges to zero uniformly in $N$  as $R\to \infty$. Consequently, combining \eqref{eq:KJ01} and \eqref{eq:KJ02}, we obtain that 
$$
\lim_{N \to \infty} \sJ_1^N =0.
$$

\item 
For $\sJ_2^N$, let $K_\eps$ be a family of molifiers in $\mR^d$ and define
$$
B_\eps(t,x) = b(t,\cdot,\rho_t(\cdot)) * K_\eps(x).
$$
Notice that $\|B_\eps\|\le \|b\|_\infty$ and for any $R>0$, $B_R:=\{x\in\mR^d\mid |x|<R \},$
\begin{align}\label{ZM01}
\lim_{\eps\to0}\|\1_{B_R}(B_\eps-b)\|_p=0.
\end{align}
Then
\begin{align*}
\sJ_2^N \leq \sJ_{21}^N(\eps) + \sJ_{22}^N(\eps) + \sJ_{23}^N(\eps),
\end{align*}
where
\begin{align*}
\sJ_{21}^N(\eps) := \mE \int_h^t  |B_\eps(s,\widetilde X^N_{\phi_N(s)}) - B_\eps(s,\widetilde X_{s}) | \dif s,
\end{align*}
\begin{align*}
\sJ_{22}^N(\eps) := \mE \int_h^t  |b(s ,\widetilde X^N_{\phi_N(s)},\rho_{s}(\widetilde X^N_{\phi_N(s)})) - B_\eps(s,\widetilde X^N_{\phi_N(s)}) | \dif s
\end{align*}
and
\begin{align*}
\sJ_{23}^N(\eps):= \mE \int_h^t  |b(s ,\widetilde X_{s},\rho_{s}(\widetilde X_{s})) - B_\eps(s,\widetilde X_{s}) | \dif s.
\end{align*} 
Thus, by \eqref{eq:AG} and results in \cite[p.\,339]{JS03}, for any $s>0$, 
$$
\widetilde X_s^N\1_{\Delta \widetilde X_s=0}(s)  \to \widetilde X_s\1_{\Delta \widetilde X_s=0}(s),\ \ \hbox{as  $N \to \infty$,  $\widetilde \mP$-a.s., }
$$
which, by the dominate convergence theorem, implies that for arbitrary fixed $\eps>0$,
\begin{align*}
\lim_{N \to \infty } \mE \int_h^t  |B_\eps(s,\widetilde X^N_{s})  - B_\eps(s,\widetilde X_{s}) | \dif s
\le &\, \mE \int_0^t\lim_{N \to \infty }  |B_\eps(s,\widetilde X^N_{s})  - B_\eps(s,\widetilde X_{s}) |\1_{ \Delta\tilde X_s=0}(s)  \dif s\\
&+2\|B_\eps\|_\infty \mE \int_0^t\1_{\Delta \tilde X_s>0}(s)\dif s  = 0,
\end{align*}
where we use the fact that for Lebsgue a.e. $s\in[0,t]$,  $\Delta\tilde X_s=0$ since $\tilde X\in\mD$.
On the other hand, by \eqref{eq:LO01}, we have
\begin{align*}
\mE \int_h^t  |B_\eps(s,\widetilde X^N_{\phi_N(s)}) - B_\eps(s,\widetilde X^N_{s}) | \dif s 
&\leq
\| \nabla B_\eps \|_\infty \int_h^t  \mE |\widetilde X^N_{\phi_N(s)} - \widetilde X^N_{s} | \dif s \\
&\lesssim  \| \nabla B_\eps \|_\infty ( |h|\|b\|_\infty + |h|^{1/\alpha}),
\end{align*}
where $h = T/N$.  Consequently, for fixed $\eps>0$,
$$
\lim_{N \to \infty}\sJ_{21}^N (\eps) =0.
$$
For $\sJ_{22}^N (\eps) $, we have
\begin{align*}
\sJ_{22}^N (\eps) 
\leq &\,
\mE \int_h^t \1_{ \{ |\widetilde X_{\phi_N(s)}^N|\leq R \} }  |b(s ,\widetilde X^N_{\phi_N(s)},\rho_{s}(\widetilde X^N_{\phi_N(s)})) - B_\eps (s ,\widetilde X^N_{\phi_N(s)}) |  \dif s\\
& +  2\|b\|_\infty \int_h^t \widetilde \mP  \( |\widetilde X_{\phi_N(s)}^N|> R \) \dif s := I_R^N(\eps) + J_R^N.
\end{align*}
Samely as \eqref{eq:ZK01}, by H\"older's inequality with $1<q < \alpha/d+1$ and $q = \frac{p}{p-1}$, we see that
\begin{align*}
I_R^N(\eps) 
\lesssim & \left[ \int_0^T 
 \int_{|y|\leq R} |b(s ,y,\rho_{s}(y) - B_\eps (s ,y) |^p \dif y   \dif s \right]^{1/p} 
\left[ \int_h^T (s-h)^{-\frac{dq}{\alpha p} } \dif s \right]^{1/q},
\end{align*}
where the implicit constant in the above $\lesssim$ is independent of $N,R$ and $\eps$. Hence, for each $R >0$, by the dominated convergence theorem and \eqref{ZM01}, we obtain
$$
\lim_{\eps \to 0} \sup_{N}I_R^N(\eps)  = 0.
$$
By \eqref{eq:KJ02}, we have $\lim_{R \to \infty}\sup_N J_R^N =0 $.  
For $\sJ_{23}^N(\eps)$, it is similar to $\sJ_{22}^N(\eps)$. 
\end{enumerate}
Combining the above calculations, we get \eqref{S3:CoB}. The proof of the existence is finished.

\medskip\noindent
{\bf(Uniqueness)}
For $i=1,2$, let $\rho_t^{(i)}$ be two densities of two {weak} solutions {$ X^{(i)}$ for DDSDE \eqref{eq:AB} respectively}:
$$
X^{(i)}_t=X_0+\int_0^t b(s,X_s^{(i)},\rho^{(i)}_s(X_s^{(i)}))\dif s+L^{(i)}_t.
$$
 For each $i$, by the well-known result (see \cite{LZ19} for example), for any $x\in\mR^d$, there is a unique weak solution $(\widetilde X^{(i)}_t(x),\widetilde L^{(i)}_t)$ for the following classical SDE with bounded drift $b(t,x,\rho^{(i)}_t(x))$ and $\alpha\in(1,2)$:
\begin{align*}
\widetilde X^{(i)}_t(x)=x+\int_0^t b(s,\widetilde X_s^{(i)}(x),\rho^{(i)}_s(\widetilde X_s^{(i)}(x)))\dif s+\widetilde L^{(i)}_t.
\end{align*}
Denote by $ p_x^{(i)}(t,\cdot)$ the density of $\widetilde X^{(i)}_t(x)$.
Then, we have
\begin{align}\label{pf:Un01}
\rho_t^{(i)}(y)=\int_{\mR^d}{ p_x^{(i)}(t,\cdot)}\rho_0(x)\dif x.
\end{align}
Moreover, by \eqref{S2:Ep} and \eqref{eq:DD04-1}, 
\begin{align}\label{pf:Un02}
{ p_x^{(i)}(t,\cdot)}\le c p_\alpha(t,x-y),\ \  \forall (t,x,y)\in(0,\infty)\times\mR^{d}\times\mR^d.
\end{align}
In the following, we prove \textbf{i)} and \textbf{ii)} respectively. 
\begin{itemize}
\item[\bf i)]
Define $u_t:=\rho^{(1)}_t-\rho^{(2)}_t$. Hence, by \eqref{eq:FG01},
\begin{align}\label{E01}
\p_t u_t=\Delta^{\alpha/2}u_t+\div(\cB(t,x)u_t), \quad u_0=0,
\end{align}
in weak sense, where 
\begin{align*}
\cB(t,x):=b(t,x,\rho^{(1)}_t(x))+\rho^{(2)}_t(x)\frac{b(t,x,\rho^{(1)}_t(x))-b(t,x,\rho^{(2)}_t(x))}{ u_t(x)},
\end{align*}
and use the convention $\frac{0}{0}=0$. In this case, we only need to prove that $u_t=0$.  By \eqref{pf:Un01} and \eqref{pf:Un02}, we have
\begin{align*}
\|\rho^{(i)}_t\|_\infty\le \|{ p_\alpha(t,\cdot)}\|_{p}\|\rho_0\|_{q}\stackrel{\eqref{eq:HG01}}{\lesssim} t^{-d/(\alpha q)}\in L^1([0,T]),\quad i=1,2,
\end{align*}
where $1/p+1/q=1$. Notice that, by the assumption \eqref{Lipb},
\begin{align}\label{ZZ00}
\|\cB(t,\cdot)\|_\infty\lesssim \|b\|_\infty+\|\rho^{(2)}_t\|_\infty \lesssim 1+ t^{-d/(\alpha q)} .
\end{align}
By Duhamel's formula and \eqref{E01}, we have
\begin{align*}
u_t(x)=\int_0^t \( p_\alpha(t-s)*\div(\cB(s)u_s)\)(x)\dif s.
\end{align*}
Hence, for $q> d/(\alpha-1)$, by \eqref{ZZ00} we have
\begin{align*}
\|u_t\|_\infty&\le \int_0^t \|\nabla p_\alpha(t-s)\|_{1}\|\cB(s)u_s\|_\infty\dif s\\
&\stackrel{\eqref{eq:HG01}}{\lesssim}\int_0^t (t-s)^{-1/\alpha}(1+s^{-d/(\alpha q)})\|u_s\|_\infty\dif s,
\end{align*}
which completes the proof by Gronwall's inequality (cf. \cite[Example 2.4]{Zh10}).
\item[\bf ii)] 
Based on the weak uniqueness result \textbf{i)} with $q=\infty$, we denote by
\begin{align*}
\rho_t:=\rho^1_t=\rho^2_t.
\end{align*}
By the well-known result (see \cite{CZZ21}*{Theorem 1.1} for example), if we have
\begin{align}\label{ZZ01}
\cA(t,x):=b(t,x,\rho_t(x))\in L^\infty([0,T];\bC^{\beta_0})
\end{align}
 for  $\beta_0\in(1-\alpha/2,1)$ and any $T>0$, then the strong uniqueness holds. Thus, it is enough to show that $\rho_t(x)\in L^\infty([0,T];\bC^{\beta_0})$ for $\beta_0>1-\alpha/2$. Unfortunately, we can not obtain it directly from $\rho_0\in \bC^{\beta_0}$ and \eqref{pf:Un01}, since $p_t(x,y)\ne p_t(x-y)$. Firstly, by \eqref{pf:Un01} and \eqref{pf:Un02}, we have
 \begin{align*}
\sup_{t\in[0,T]}\|\rho_t\|_\infty \lesssim \|\rho_0\|_\infty\sup_{t\in[0,T]}\int_{\mR^d} p_\alpha(t,y)\dif y=\|\rho_0\|_\infty,
\end{align*}
 which implies that
 \begin{align*}
 \bar{b}(t,x):=b(t,x,\rho_t(x))\rho_t(x)\in L^\infty([0,T]\times \mR^d).
\end{align*}
Hence, by \eqref{eq:FG01}, \eqref{S2:Sch}, \eqref{S2:Bern} and \eqref{LJ1}, for any $T>0$, there is a constant $c_T$ such that for all $t\in[0,T]$,
\begin{align}\label{eq:QO01}
\|\rho_t\|_{\bC^{(\alpha-1)\wedge\beta_0}}
\le c_T\(\|\rho_0\|_{\bC^{\beta_0}}+\|\bar{b}\|_{\mL^\infty_T}\).
\end{align}
If $\alpha-1>\beta_0$, \eqref{ZZ01} is straightforward. Otherwise,
by \eqref{Lipb} and \eqref{eq:QO01}, we have
\begin{align*}
\|\bar{b}\|_{L^\infty_T(\bC^{\alpha-1})}<\infty.
\end{align*}
Thus, by \eqref{eq:FG01}, \eqref{S2:Sch}, \eqref{S2:Bern} and \eqref{LJ1} again, there is a constant $c^{(2)}_T$ such that for all $t\in[0,T]$,
\begin{align*}
\|\rho_t\|_{\bC^{(2(\alpha-1))\wedge \beta_0}}\le c^{(2)}_T.
\end{align*}
By induction, there are a $N\in\mN$ with
\begin{align*}
(\alpha-1)N>\beta_0
\end{align*}
and a constant $c^{(N)}_T$ such that
\begin{align*}
\|\rho_t\|_{\bC^{\beta_0}}\le c^{(N)}_T,\ \ \forall t\in[0,T].
\end{align*}
Then, we have \eqref{ZZ01} and complete the proof.
\end{itemize}
\end{proof}

\section*{}
\subsection*{Acknowledgments}
We are deeply grateful to Prof. Xicheng Zhang and Prof. Fuke Wu for their valuable suggestions and correcting some errors. 

%%%%%%%%%%%%%%%%%%
%                 APPENDIX                 %
%%%%%%%%%%%%%%%%%%

%\newpage
%\begin{appendix}
%\section{Appendix}
%\Appendix
%\end{appendix}

%\bibliographystyle{BIBbase/amsxport}
%\bibliography{BIBbase/referenceWu.bib}


\end{document}